\documentclass[12pt]{amsart}
\usepackage{amssymb}

\usepackage[]{amsfonts}

\def\underset#1#2{\mathrel{\mathop{\kern0pt #2}\limits_{#1}}}

\def\overset#1#2{\mathrel{\mathop{\kern0pt #2}\limits^{#1}}}

\def\couleur(#1 #2 #3)
	{
\special{ps::[end]}}

\def\sqr#1#2{{\vcenter{\vbox{\hrule height.#2pt
             \hbox{\vrule width.#2pt height#1pt \kern#1pt
             \vrule width.#2pt}
             \hrule height.#2pt}}}}

\def\st{\mathinner{\mkern1mu\raise1pt\hbox{.}				
		   \mkern1mu\raise4pt\hbox{.}
		   \mkern1mu\raise1pt\hbox{.}
		 }
         }

\def\bx#1{\setbox1=\hbox{\kern3pt{#1}\kern3pt}				
 \dimen1=\ht1 \advance\dimen1 by 3pt \dimen2=\dp1 \advance\dimen2 by 3pt
 \setbox1=\hbox{\vrule height\dimen1 depth\dimen2\box1\vrule}%
 \setbox1=\vbox{\hrule\box1\hrule}%
 \advance\dimen1 by .4pt \ht1=\dimen1
 \advance\dimen2 by .4pt \dp1=\dimen2 \box1\relax}

\def\k#1{\kern#1em}
\def\vci{\vrule  width.02em height1.47ex depth-.0ex}				
\def\11{{\rm\k{.2}\vci\k{-.37}1}}

\def\liminf{\mathop{\underline{\rm lim}}}
\def\limsup{\mathop{\overline{\rm lim}}}

\newtheorem{Theorem}{Theorem}[section]
\newtheorem{Lemma}[Theorem]{Lemma}
\newtheorem{Proposition}[Theorem]{Proposition}
\newtheorem{Remark}[Theorem]{Remark}

\begin{document}
\title{Numerical solution of a parabolic system with blow-up of the solution}
\author{Marie-Noelle LeRoux}
\address{UNIVERSITE BORDEAUX1, Institut de Math{\'e}matiques de Bordeaux, 
UMR 5251, 351, Cours de la Lib{\'e}ration,33405, TALENCE Cedex}
\email{Marie-Noelle.Leroux@math.u-bordeaux1.fr}
\maketitle
\begin{abstract} {
In this paper, the author proposes a numerical method to solve a parabolic 
system of two quasilinear equations of nonlinear heat conduction with 
sources. The solution of this system may blow up in finite time. It is 
proved that the numerical solution also may blow up in finite time and 
an estimate of this time is obtained. The convergence of the scheme is 
obtained for particular values of the parameters. \ \par
}\end{abstract}
\ \par
\section{Introduction}
\setcounter{equation}{0}\ \par
{\hskip 1.8em}The purpose of this paper is to study the numerical behavior 
of the solution of a nonlinear reaction diffusion system with nonlinear 
source terms.\ \par
\ \par
\ \par
{\hskip 1.8em}Let $\Omega $ a smooth bounded domain in ${\mathbb{R}}^{d}$. 
We consider the system:\ \par

\begin{displaymath} 
\left\lbrace{\begin{matrix}{u_{1t}-\Delta u_{1}^{\nu +1}=\alpha v_{1}^{\mu 
+1},\ x\in \Omega ,\ t>0}\cr {v_{1t}-\Delta v_{1}^{\mu +1}=\alpha u_{1}^{\nu 
+1},\ x\in \Omega ,\ t>0}\cr {u_{1}(x,0)=u_{10}(x)>0,\ x\in \Omega }\cr 
{v_{1}(x,0)=v_{10}(x)>0,\ x\in \Omega }\cr {u_{1}=v_{1}=0,\ x\in \partial 
\Omega ,t>0}\cr \end{matrix}}\right. \end{displaymath} \ \par
with $\displaystyle \nu ,\mu >0,\ \alpha \geq 0$. \ \par
Samarskii and al. ~\cite{Samarskii1} have studied this system and obtained 
the following results:\ \par
If $\lambda _{1}\ $denotes the first eigenvalue of the Dirichlet problem: 
$-\Delta \rho =\lambda \rho ,\ x\in \Omega ,\ \rho =0,x\in \partial \Omega 
$ and if $\alpha >\lambda _{1}$, the problem has no global solutions 
and there exists $T_{0}>0$ such that \ \par

\begin{displaymath} 
\underset{t{\longrightarrow}T_{0}-}{\limsup }\left({\left\Vert{u_{1}^{\nu 
+1}(t,.)}\right\Vert ^{2}_{L^{2}(\Omega )}+\left\Vert{v_{1}^{\mu +1}(t,.)}\right\Vert 
^{2}_{L^{2}(\Omega )}}\right) =+\infty \end{displaymath} \ \par
In ~\cite{mnLeRoux1}, ~\cite{mnLeRoux3}, we have proposed a numerical 
method to solve a quasilinear parabolic equation with blow-up of the 
solution. The numerical solution is computed by using the function $u=u_{1}^{\nu +1};\ $so 
the nonlinearity is reported on the derivative in time. This solution 
has the same properties as the exact solution, in particular blow-up 
in finite time. We generalize  this method to the system of two equations.\ 
\par
\ \par
\ \par
{\hskip 1.8em}For what follows, it is more convenient to work with a 
transformed equation. Let $\displaystyle u=u_{1}^{\nu +1}$, $\displaystyle 
v=v_{1}^{\mu +1}$, $m=\frac{1}{\nu +1},\ p=\frac{1}{\mu +1}$; then we 
get $m,p\in ]0,1[$ and we suppose that $p\leq m$ (or $\mu \geq \nu $).\ 
\par
Then $(u,v)$ satisfies the following system: \ \par

\begin{equation} 
\displaystyle \left\lbrace{\displaystyle \begin{matrix}{\displaystyle 
mu^{m-1}u_{t}+Au=\alpha v}\cr {\displaystyle pv^{p-1}v_{t}+Av=\alpha 
u}\cr {\displaystyle u(x,0)=u_{10}^{\nu +1}(x)=u_{0}(x),\ x\in \Omega 
}\cr {\displaystyle v(x,0)=v_{10}^{\mu +1}(x)=v_{0}(x),\ x\in \Omega 
}\cr \end{matrix}}\right. \label{systeme0}
\end{equation} \ \par
where $A$ is the operator $-\Delta $ of domain $D(A)=H_{0}^{1}(\Omega )\cap H^{2}(\Omega )$.\ 
\par
\ \par
An outline of the paper is as follows: In Section 2, we study the asymptotic 
behavior of the solution. In Section 3, we define a numerical scheme 
and prove the existence of the solution of this scheme. The section 4 
is devoted to the properties of the scheme, in particular, the existence 
of a numerical blow-up time in the case $\alpha >\lambda _{1}.$ Finally, 
in Section 5, we study the particular case $p=m$ and prove the convergence 
of the scheme in that case for a class of initial conditions.\ \par
\ \par
\section{Asymtotic behavior of the solution}
\setcounter{equation}{0}\ \par
Given $u_{0},v_{0}\in L^{\infty }(\Omega )$, a couple $(u,v)$ is a weak 
solution of ~(\ref{systeme0}) on $[0,T]$ if $u,v\in L^{\infty }((0,T)\times \Omega )$ 
and 
\begin{displaymath} 
\int_{0}^{t}{\int_{\Omega }^{}{(u^{m}\phi _{t}-uA\phi +\alpha v\phi )dx}}dt=\int_{\Omega 
}^{}{u(x,t)\phi (x)dx}-\int_{\Omega }^{}{u_{0}(x)\phi (x)dx}\end{displaymath} 
\ \par

\begin{displaymath} 
\int_{0}^{t}{\int_{\Omega }^{}{(v^{p}\phi _{t}-vA\phi +\alpha u\phi )dx}}dt=\int_{\Omega 
}^{}{v(x,t)\phi (x)dx}-\int_{\Omega }^{}{v_{0}(x)\phi (x)dx}\end{displaymath} 
\ \par
for all $\phi \in C^{2}((0,T)\times \Omega )\cap C^{1}([0,T]\times \overline{\Omega 
}),\phi (x,t)=0\ for\ x\in \partial \Omega $.  \ \par
\ \par
{\hskip 1.8em}This problem admits a local solution  and from the maximum 
principle, we get $u(t),v(t)>0\ $in $\Omega ,\ 0<t<T$.\ \par
We prove the following results: \ \par
$\iota $) if $\alpha >\lambda _{1}$, the solution blows up in finite 
time$\ T$; we get: \ \par
\ \par

\begin{displaymath} 
\underset{\displaystyle t{\longrightarrow}T_{-}}{\displaystyle \lim }\displaystyle 
\int_{\Omega }^{}{\displaystyle \left({\displaystyle \frac{\displaystyle 
m}{\displaystyle m+1}u^{m+1}(t)+\frac{\displaystyle p}{\displaystyle 
p+1}v^{p+1}(t)}\right) dx}=+\infty \end{displaymath} \ \par
\ \par
$\iota \iota $) if $\alpha <\lambda _{1}$, the problem has a global solution 
which tends to $0$ when $t{\longrightarrow}\infty $ \ \par
\ \par
$\iota \iota \iota $) if $\alpha =\lambda _{1}$, the problem has a global 
solution $(u,v)\ $which tends to $\theta \rho _{1}$ when $t{\longrightarrow}\infty $ 
where $\rho _{1}$ is the first eigenfunction of $A$ ($A\rho _{1}=\lambda _{1}\rho _{1}$ 
and $\left\Vert{\rho _{1}}\right\Vert _{L^{1}(\Omega )}=1$) and $\theta $ 
is a constant depending on the initial condition.\ \par
\ \par
We introduce the functions $\Phi $, $Z$ defined on $[0,T]$ by \ \par

\begin{displaymath} 
\Phi (t)=\int_{\Omega }^{}{\left({\frac{m}{m+1}u^{m+1}(t)+\frac{p}{p+1}v^{p+1}(t)}\right) 
dx}\end{displaymath} \ \par

\begin{equation} 
Z(t)=\ \Phi (t)^{\frac{m-1}{m+1}}=\left({\int_{\Omega }^{}{\left({\frac{m}{m+1}u^{m+1}(t)+\frac{p}{p+1}v^{p+1}(t)}\right) 
dx}}\right) ^{\frac{m-1}{m+1}}\label{systeme1}
\end{equation} \ \par
and the functional defined on $H_{0}^{1}(\Omega )\times H_{0}^{1}(\Omega )$by 
\ \par
\ \par

\begin{equation} 
J(u,v)=\int_{\Omega }^{}{(\left\vert{\nabla u}\right\vert ^{2}+\left\vert{\nabla 
v}\right\vert ^{2}-2\alpha uv)dx}\label{systeme2}
\end{equation} \ \par
\ \par
\begin{Lemma} The function $\Phi $ is convex and the function $Z$ is 
concave\label{systeme6}\ \par
\end{Lemma}
\textsl{Proof:} We prove that the second derivative of $\Phi $ is nonnegative.\ 
\par
We have: $\Phi '(t)=\int_{\Omega }^{}{(mu^{m}u_{t}+pv^{p}v_{t})dx}$\ 
\par
By multiplying the first equation of ~(\ref{systeme0}) by $u$,  the second 
by $v$ and integrating over $\Omega $, we get: 
\begin{displaymath} 
\Phi '(t)=\displaystyle \int_{\Omega }^{}{(mu^{m}u_{t}+pv^{p}v_{t})dx}=-\displaystyle 
\int_{\Omega }^{}{((Au-\alpha v)u+(Av-\alpha u)v)dx}=\end{displaymath} 
 
\begin{displaymath} 
-\int_{\Omega }^{}{(\left\vert{\nabla u}\right\vert ^{2}+\left\vert{\nabla 
v}\right\vert ^{2}-2\alpha uv)dx}=-J(u,v).\end{displaymath} \ \par
Then we deduce:  
\begin{displaymath} 
\Phi ''(t)=-2\int_{\Omega }^{}{((Au-\alpha v)u_{t}+(Av-\alpha u)v_{t})dx}\end{displaymath} 
\ \par

\begin{displaymath} 
=2\int_{\Omega }^{}{\left({\frac{1}{m}(Au-\alpha v)^{2}u^{1-m}+\frac{1}{p}(Av-\alpha 
u)^{2}v^{1-p}}\right) dx}\end{displaymath} \ \par
and $\Phi ''(t)\geq 0,\ \forall t\geq 0$. The function $\Phi $ is convex.\ 
\par
\ \par
Besides, we have:
\begin{displaymath} 
Z'(t)=\frac{m-1}{m+1}\Phi (t)^{-2/(m+1)}\Phi '(t)\end{displaymath} \ 
\par
and $Z''(t)=\frac{1-m}{m+1}\left({\Phi (t)}\right) ^{-2/(m+1)-1}\left({\frac{2}{m+1}(\Phi 
'(t))^{2}-\Phi (t)\Phi ''(t)}\right) .$\ \par
\ \par
By using Cauchy-Schwarz inequality, we get:
\begin{displaymath} 
\int_{\Omega }^{}{(Au-\alpha v)udx}\leq \left({\frac{m+1}{m}\int_{\Omega 
}^{}{(Au-\alpha v)^{2}u^{1-m}dx}}\right) ^{1/2}\left({\int_{\Omega }^{}{\frac{m}{m+1}u^{m+1}dx}}\right) 
^{1/2}\end{displaymath} \ \par
and
\begin{displaymath} 
\int_{\Omega }^{}{(Av-\alpha u)vdx}\leq \left({\frac{p+1}{p}\int_{\Omega 
}^{}{(Av-\alpha u)^{2}v^{1-p}dx}}\right) ^{1/2}\left({\int_{\Omega }^{}{\frac{p}{p+1}v^{p+1}dx}}\right) 
^{1/2}\end{displaymath} \ \par
We deduce:
\begin{displaymath} 
\ \end{displaymath}  
\begin{displaymath} 
(\Phi '(t))^{2}\leq \displaystyle \left({\displaystyle \frac{\displaystyle 
m+1}{\displaystyle m}\displaystyle \int_{\Omega }^{}{(Au-\alpha v)^{2}u^{1-m}dx}+\frac{\displaystyle 
p+1}{\displaystyle p}\displaystyle \int_{\Omega }^{}{(Av-\alpha u)^{2}v^{1-p}dx}}\right) 
\Phi (t)\end{displaymath} \ \par
and 
\begin{displaymath} 
Z''(t)\leq 2\frac{\displaystyle (1-m)}{\displaystyle (m+1)^{2}}\ \frac{\displaystyle 
p-m}{\displaystyle p}\displaystyle \left({\displaystyle \displaystyle 
\int_{\Omega }^{}{\displaystyle \left({\displaystyle \frac{\displaystyle 
m}{\displaystyle m+1}u^{m+1}+\frac{\displaystyle p}{\displaystyle p+1}v^{p+1}}\right) 
dx}}\right) ^{-2/(m+1)}\displaystyle \int_{\Omega }^{}{(Av-\alpha u)^{2}v^{1-p}dx}.\end{displaymath} 
\ \par
Since $p\leq m$, the function $Z''$ is nonpositive , that is the function 
$Z$ is concave.\ \par
\ \par
In the two next propositions, we prove that in the case $\alpha >\lambda _{1}$, 
the solution $(u,v)${\it{ }}of ~(\ref{systeme0}){\it{ }}blows up in finite 
time and obtain estimates of this time for a class of initial conditions.\ 
\par
\ \par
\begin{Proposition} If $\alpha >\lambda _{1}$ and if the initial condition 
satisfies $J(u_{0},v_{0})<0$, the solution blows up in finite time $T\ $such 
that \ \par

\begin{equation} 
T<\frac{\displaystyle 1+m}{\displaystyle 1-m}\frac{\displaystyle \displaystyle 
\int_{\Omega }^{}{\displaystyle \left({\displaystyle \frac{\displaystyle 
m}{\displaystyle m+1}u_{0}^{m+1}+\frac{\displaystyle p}{\displaystyle 
p+1}v_{0}^{p+1}}\right) dx}}{\displaystyle -\displaystyle \int_{\Omega 
}^{}{\displaystyle \left({\displaystyle \displaystyle \left\vert{\displaystyle 
\nabla u_{0}}\right\vert ^{2}+\displaystyle \left\vert{\displaystyle 
\nabla v_{0}}\right\vert ^{2}-2\alpha u_{0}v_{0}}\right) dx}}.\label{systeme8}
\end{equation} \ \par
 \ \par
\end{Proposition}
\textsl{Proof: }The function $Z'$ is nonincreasing, so we get: \ \par

\begin{equation} 
Z(t)-Z(0)\leq tZ'(0)\label{systeme9}
\end{equation} \ \par
and $\displaystyle Z'(0)=\frac{\displaystyle 1-m}{\displaystyle 1+m}\ (\Phi (0))^{\displaystyle 
\frac{\displaystyle -2}{\displaystyle m+1}}J(u_{0},v_{0})$.\ \par
\ \par
The inequality ~(\ref{systeme9}) may be written as\ \par

\begin{displaymath} 
\left({\int_{\Omega }^{}{\left({\frac{m}{m+1}u(t)^{m+1}+\frac{p}{p+1}v(t)^{p+1}}\right) 
dx}}\right) ^{\frac{m-1}{m+1}}\leq \left({\int_{\Omega }^{}{\left({\frac{m}{m+1}u_{0}^{m+1}+\frac{p}{p+1}v_{0}^{p+1}}\right) 
dx}}\right) ^{\frac{m-1}{m+1}}\end{displaymath} \ \par

\begin{equation} 
\left({1+\frac{1-m}{m+1}t\left({\int_{\Omega }^{}{\left({\frac{m}{m+1}u_{0}^{m+1}+\frac{p}{p+1}v_{0}^{p+1}}\right) 
dx}}\right) ^{-1}J(u_{0},v_{0})}\right) .\label{systeme3}
\end{equation} \ \par
\ \par
{\hskip 1.8em}If $\alpha >\lambda _{1}$ , the set $S=\{(u_{0},v_{0})\in H_{0}^{1}(\Omega )/\ J(u_{0},v_{0})<0\}$ 
is not empty.  If $J(u_{0},v_{0})<0$, the right member of ~(\ref{systeme3}) 
becomes equal to zero for a finite time and so the solution blows up 
in a finite time $T$ such that \ \par

\begin{displaymath} 
T\leq \frac{\displaystyle m+1}{\displaystyle 1-m}\ \frac{\displaystyle 
\displaystyle \int_{\Omega }^{}{\displaystyle \left({\displaystyle \frac{\displaystyle 
m}{\displaystyle m+1}u_{0}^{m+1}+\frac{\displaystyle p}{\displaystyle 
p+1}v_{0}^{p+1}}\right) dx}}{\displaystyle -J(u_{0},v_{0})}.\end{displaymath} 
\ \par
\begin{Proposition} : If $\alpha >\lambda _{1}$, the solution blows up 
in finite time $T$ and satisfies the inequality:\ \par

\begin{displaymath} 
\left({\int_{\Omega }^{}{\left({\frac{m}{m+1}u^{m+1}(t)+\frac{p}{p+1}v^{p+1}(t)}\right) 
dx}}\right) ^{\frac{1}{m+1}}\leq \left({\frac{T}{T-t}}\right) ^{^{\frac{1}{1-m}}}\left({\int_{\Omega 
}^{}{\left({\frac{m}{m+1}u_{0}^{m+1}+\frac{p}{p+1}v_{0}^{p+1}}\right) 
dx}}\right) ^{\frac{1}{m+1}}.\label{systeme5}\end{displaymath} \ \par
\end{Proposition}
\textsl{Proof:} If $J(u_{0},v_{0})<0,$ from Proposition 2.3,  we know 
that the solution blows up in finite time. If $J(u_{0},v_{0})\geq 0$, 
the result is obtained by using the same proof as Friedman-McLeod in 
~\cite{Friedman}.\ \par
\ \par
{\hskip 1.8em}Since the function $Z$  is concave, it satisfies \ \par
$\displaystyle Z(t)\geq \frac{\displaystyle s-t}{\displaystyle s}Z(0)+\frac{\displaystyle 
t}{\displaystyle s}Z(s),\ 0\leq t\leq s<T$  and $\displaystyle \underset{\displaystyle s{\longrightarrow}T_{-}}{\displaystyle \lim }Z(s)=0$.\ 
\par
We deduce 
\begin{displaymath} 
\left({\int_{\Omega }^{}{\left({\frac{m}{m+1}u^{m+1}(t)+\frac{p}{p+1}v^{p+1}(t)}\right) 
dx}}\right) ^{\frac{m-1}{m+1}}\geq \frac{T-t}{T}\left({\int_{\Omega }^{}{\left({\frac{m}{m+1}u_{0}^{m+1}+\frac{p}{p+1}v_{0}^{p+1}}\right) 
dx}}\right) ^{\frac{m-1}{m+1}},\end{displaymath} \ \par
that is  the inequality ~(\ref{systeme5}).\ \par
\ \par
\begin{Proposition} If $\alpha <\lambda _{1}$, the solution tends to 
$0$ when$\ t{\longrightarrow}+\infty .$\ \par
\end{Proposition}
\ \par
\textsl{Proof: } If $\lambda _{1}\geq \alpha ,$ we get: $J(u,v)\geq (\lambda _{1}-\alpha )\int_{\Omega }^{}{(u^{2}+v^{2})dx}\geq 
0\ $ and $\Phi '(t)\leq 0$. Since the function $\Phi $ is convex, the 
function $\Phi '$ is nondecreasing and $\underset{t{\longrightarrow}+\infty }{\lim }\Phi '(t)=l\leq 0$.\ 
\par
We deduce: $\Phi (t)\leq \Phi (0)+lt,\ \forall t>0$; since $\Phi (t)\geq 0,$ 
we obtain that $l=0$ and \ \par

\begin{equation} 
\underset{t{\longrightarrow}+\infty }{\lim }J(u(t),v(t))=0.\label{systeme15}
\end{equation} \ \par
If $\lambda _{1}-\alpha >0,\ $ then, we get  $\underset{t{\longrightarrow}+\infty }{\lim }\int_{\Omega }^{}{(u^{2}+v^{2})dx}=0$ 
and the solution tends to $0$ when $t{\longrightarrow}+\infty .$\ \par
\ \par
\begin{Proposition} If $\lambda _{1}=\alpha $, then $\underset{t{\longrightarrow}+\infty }{\lim }u(t)=\underset{t{\longrightarrow}+\infty 
}{\lim }v(t)=\theta \rho _{1}$ in $L^{2}(\Omega )$ where $\theta $ is 
a constant depending on the initial conditions.\label{systeme33}\ \par
\end{Proposition}
\ \par
\textsl{Proof:}  If $\lambda _{1}=\alpha $, we get from ~(\ref{systeme15}) 
 that $\underset{t{\longrightarrow}+\infty }{\lim }J(u(t),v(t))=0$ and 
$\Phi (t)$ is bounded for $t\geq 0.$\ \par
By interpolation, we obtain\ \par
$\displaystyle \displaystyle \int_{\Omega }^{}{u^{2}dx}\leq \delta \displaystyle \int_{\Omega 
}^{}{\displaystyle \left\vert{\displaystyle \nabla u}\right\vert ^{2}dx}+C_{1}(\delta 
)\displaystyle \left({\displaystyle \displaystyle \int_{\Omega }^{}{u^{m+1}dx}}\right) 
^{2/(m+1)},$\ \par
\ \par
$\displaystyle \displaystyle \int_{\Omega }^{}{v^{2}dx}\leq \delta \displaystyle \int_{\Omega 
}^{}{\displaystyle \left\vert{\displaystyle \nabla v}\right\vert ^{2}dx}+C_{2}(\delta 
)\displaystyle \left({\displaystyle \displaystyle \int_{\Omega }^{}{v^{p+1}dx}}\right) 
^{2/(p+1)},$\ \par
where $C_{1}(\delta )\ $and $C_{2}(\delta )$ are constants depending 
on $\Omega $ and $m$ and $p$ respectively and we get:\ \par
\ \par
$\displaystyle \displaystyle \int_{\Omega }^{}{(u^{2}+v^{2})dx}\leq \delta J(u(t),v(t))+2\delta 
\lambda _{1}\displaystyle \int_{\Omega }^{}{uvdx}$ $\displaystyle +C_{1}(\delta )\displaystyle \left({\displaystyle \displaystyle \int_{\Omega 
}^{}{u^{m+1}dx}}\right) ^{2/(m+1)}+C_{2}(\delta )\displaystyle \left({\displaystyle 
\displaystyle \int_{\Omega }^{}{v^{p+1}dx}}\right) ^{2/(p+1)}$.\ \par
\ \par
We deduce \ \par
$\displaystyle (1-\delta \lambda _{1})\displaystyle \int_{\Omega }^{}{(u^{2}+v^{2})dx}\leq 
\delta J(u(t),v(t))+C_{1}(\delta )\displaystyle \left({\displaystyle 
\displaystyle \int_{\Omega }^{}{u^{m+1}dx}}\right) ^{2/(m+1)}+C_{2}(\delta 
)\displaystyle \left({\displaystyle \displaystyle \int_{\Omega }^{}{v^{p+1}dx}}\right) 
^{2/(p+1)}.$\ \par
\ \par
Since $J(u(t),v(t))$ is bounded, if $\delta $ is chosen such that $1-\delta \lambda _{1}>0$, 
we deduce that $u(t)$ and $v(t)$ are uniformly bounded in $L^{2}(\Omega )$ 
and in $H_{0}^{1}(\Omega )$. So, we can extract subsequences $t_{n}{\longrightarrow}+\infty $ 
such that $u(t_{n})$ and $v(t_{n})$ converge weakly in $H_{0}^{1}(\Omega )$ 
and strongly in $L^{2}(\Omega )$ to $z_{1}$ and $z_{2}$ respectively.\ 
\par
\ \par
We have : $J(z_{1},z_{2})\leq \ \liminf J(u(t),v(t))=0$, that is $\int_{\Omega }^{}{\left({\left\vert{\nabla z_{1}}\right\vert ^{2}+\left\vert{\nabla 
z_{2}}\right\vert ^{2}-2\lambda _{1}z_{1}z_{2}}\right) dx}=0\ $or\ \par

\begin{displaymath} 
\int_{\Omega }^{}{\left({\left\vert{\nabla z_{1}}\right\vert ^{2}-\lambda 
_{1}z_{1}^{2}}\right) dx}+\int_{\Omega }^{}{\left({\left\vert{\nabla 
z_{2}}\right\vert ^{2}-\lambda _{1}z_{2}^{2}}\right) dx}+\lambda _{1}\int_{\Omega 
}^{}{\left({z_{1}-z_{2}}\right) ^{2}dx}=0\end{displaymath} \ \par
We deduce that $z_{1}=z_{2}=\theta \rho _{1}.$ \ \par
\ \par
By multiplying the two equations of ~(\ref{systeme0}) by $\rho _{1}$ 
and integrating over $\Omega $, we get:\ \par
$\displaystyle \frac{\displaystyle d}{\displaystyle dt}\displaystyle \left({\displaystyle 
\displaystyle \int_{\Omega }^{}{(u^{m}+v^{p}})\rho _{1}dx}\right) =0$ 
and then\ \par
$\displaystyle \displaystyle \int_{\Omega }^{}{}(u^{m}(t)+v^{p}(t))\rho _{1}dx=\displaystyle 
\int_{\Omega }^{}{(u_{0}^{m}+v_{0}^{m})\rho _{1}dx},\ $for all $t>0$.\ 
\par
\ \par
Hence we obtain\ \par

\begin{displaymath} 
\int_{\Omega }^{}{(\theta ^{m}\rho _{1}^{m+1}+\theta ^{p}\rho _{1}^{p+1})dx}=\int_{\Omega 
}^{}{(u_{0}^{m}+v_{0}^{p})\rho _{1}dx}.\end{displaymath} \ \par
and there exists a unique positive value of $\theta $ satisfying this 
equation; we deduce the proposition\ \par
\ \par
\section{Definition of a numerical scheme}
\setcounter{equation}{0}\ \par
The classical Euler scheme cannot blow up in a finite time, so we generalize 
here to a system the numerical scheme used in ~\cite{mnLeRoux1}.\ \par
The first equation of ~(\ref{systeme0}) may be written:\ \par
\ \par

\begin{displaymath} 
-\frac{m}{1-p}u^{m-p+1}\frac{d}{dt}(u^{p-1})+Au=\alpha v\end{displaymath} 
\ \par
and the second: 
\begin{displaymath} 
-\frac{p}{1-p}v\frac{d}{dt}(v^{p-1})+Av=\alpha u.\end{displaymath} \ 
\par
\ \par
So, we discretize the two derivatives in time in the same manner:\ \par
If $(u_{n},v_{n})$ is the approximate solution at the time level $t_{n}=n\Delta t$, 
( where  $\Delta t$ is the time step), the approximate solution at the 
time level $t_{n+1}$ is solution of the system:\ \par

\begin{equation} 
\frac{m}{1-p}u_{n+1}u_{n}^{m-p}(u_{n}^{p-1}-u_{n+1}^{p-1})+\Delta tAu_{n+1}=\alpha 
\Delta tv_{n+1},\label{systeme11}
\end{equation} \ \par

\begin{equation} 
\frac{p}{1-p}v_{n+1}(v_{n}^{p-1}-v_{n+1}^{p-1})+\Delta tAv_{n+1}=\alpha 
\Delta tu_{n+1}.\label{systeme12}
\end{equation} \ \par
We prove first that the system ~(\ref{systeme11}), ~(\ref{systeme12}) 
has a unique positive solution if $u_{n},v_{n}$ are positive in $\Omega $.\ 
\par
We need several lemmas. For what follows, we denote: $\left\Vert{v}\right\Vert _{r}=\left\Vert{v}\right\Vert _{L^{r}(\Omega 
)}.$\ \par
\begin{Lemma} If the functions $u_{n}$ and $v_{n}$ are positive on $\Omega ,\ $continuous 
in $\overline{\Omega }$ and satisfy the condition: \ \par

\begin{equation} 
\left\Vert{u_{n}}\right\Vert _{\infty }^{1-m}\left\Vert{v_{n}}\right\Vert 
_{\infty }^{1-p}<\frac{mp}{\alpha ^{2}(1-p)^{2}\Delta t^{2}}\label{systeme13}
\end{equation} \ \par
then the system ~(\ref{systeme11}), ~(\ref{systeme12}) has a positive 
solution $u_{n+1},v_{n+1}\in C^{2}(\overline{\Omega })$ \label{systeme17}\ 
\par
\end{Lemma}
\ \par
\textsl{Proof: } Consider the functional defined on $H_{0}^{1}(\Omega )\times H_{0}^{1}(\Omega )$ 
by: \ \par

\begin{equation} 
J_{n}(u,v)=\int_{\Omega }^{}{(\left\vert{\nabla u}\right\vert ^{2}+\left\vert{\nabla 
v}\right\vert ^{2}-2\alpha uv)dx}+\frac{1}{(1-p)\Delta t}\int_{\Omega 
}^{}{(mu_{n}^{m-1}u^{2}+pv_{n}^{p-1}v^{2})dx}\label{systeme24}
\end{equation} \ \par
and let us denote $\displaystyle K=\displaystyle \left\lbrace{\displaystyle \displaystyle \left.{\displaystyle 
(u,v)\in (H_{0}^{1}(\Omega ))^{2}/\displaystyle \int_{\Omega }^{}{(mu_{n}^{m-p}\displaystyle 
\left\vert{\displaystyle u}\right\vert ^{p+1}+p\displaystyle \left\vert{\displaystyle 
v}\right\vert ^{p+1})dx=1}}\right\rbrace }\right. .\ \ $\ \par
\ \par
Since $(u_{n},v_{n})$ satisfies ~(\ref{systeme13}), we get : $mu_{n}^{m-1}u^{2}+pv_{n}^{p-1}v^{2}-2\alpha (1-p)\Delta tuv\geq 0$ 
and\ \par
$J_{n}(u,v)\geq 0,\ J_{n}(\left\vert{u}\right\vert ,\left\vert{v}\right\vert 
)\leq J_{n}(u,v),\ \forall u,v\in H_{0}^{1}(\Omega ).$\ \par
Now, we consider the problem:   
\begin{equation} 
\underset{(u,v)\in K}{\min }J_{n}(u,v)\label{systeme14}
\end{equation} \ \par
From the preceding remark, the solution of ~(\ref{systeme14}), if it 
exists, will be positive.\ \par
We denote $\phi (u,v)=\int_{\Omega }^{}{(mu_{n}^{m-p}u^{p+1}+pv^{p+1})dx}.$\ 
\par
If $({\hat u},\hat v)\ $is a solution of problem ~(\ref{systeme14}), 
there exists $\lambda \in {\mathbb{R}}$ such that: \ \par

\begin{displaymath} 
\displaystyle \left\lbrace{\displaystyle \begin{matrix}{\displaystyle 
\frac{\displaystyle \partial J_{n}({\hat u},\hat v)}{\displaystyle \partial 
u}-\lambda \frac{\displaystyle \partial \phi }{\displaystyle \partial 
u}({\hat u},\hat v)=0}\cr {\displaystyle \frac{\displaystyle \partial 
J_{n}({\hat u},\hat v)}{\displaystyle \partial v}-\lambda \frac{\displaystyle 
\partial \phi }{\displaystyle \partial v}({\hat u},\hat v)=0}\cr {\displaystyle 
\phi ({\hat u},\hat v)=1}\cr \end{matrix}}\right. \end{displaymath} \ 
\par
Then, we get for any $\psi \in H_{0}^{1}(\Omega )$,\ \par

\begin{displaymath} 
\int_{\Omega }^{}{(\nabla {\hat u}\ \nabla \psi -\alpha \hat v\psi )dx}+\frac{m}{(1-p)\Delta 
t}\int_{\Omega }^{}{u_{n}^{m-1}{\hat u}\psi dx}=\lambda m(p+1)\int_{\Omega 
}^{}{u_{n}^{m-p}{\hat u}^{p}\psi dx},\end{displaymath} \ \par

\begin{displaymath} 
\int_{\Omega }^{}{(\nabla \hat v\ \nabla \psi -\alpha {\hat u}\psi )dx}+\frac{p}{(1-p)\Delta 
t}\int_{\Omega }^{}{v_{n}^{p-1}\hat v\psi dx}=\lambda p(p+1)\int_{\Omega 
}^{}{\hat v^{p}\psi dx}\end{displaymath} \ \par
and
\begin{displaymath} 
\int_{\Omega }^{}{(mu_{n}^{p-1}{\hat u}^{p+1}+p\hat v^{p+1})dx}=1.\end{displaymath} 
\ \par
Hence, $({\hat u},\hat v)$ satisfies the equalities:\ \par

\begin{displaymath} 
A{\hat u}-\alpha \hat v+\frac{m}{(1-p)\Delta t}u_{n}^{m-1}{\hat u}=\lambda 
m(p+1)u_{n}^{m-p}{\hat u}^{p},\end{displaymath} \ \par

\begin{displaymath} 
A\hat v-\alpha {\hat u}+\frac{p}{(1-p)\Delta t}v_{n}^{p-1}\hat v=\lambda 
p(p+1)\hat v^{p}\end{displaymath}  and $({\hat u},\hat v)\in C^{2}(\Omega ).$ 
Besides $J_{n}({\hat u},\hat v)=\lambda (p+1).$\ \par
A solution of ~(\ref{systeme11}), ~(\ref{systeme12}) is then defined 
by $u_{n+1}=\gamma {\hat u},\ v_{n+1}=\gamma \hat v;$ we get
\begin{displaymath} 
Au_{n+1}-\alpha v_{n+1}+\frac{m}{(1-p)\Delta t}u_{n}^{m-1}u_{n+1}=\lambda 
\gamma ^{1-p}m(p+1)u_{n}^{m-p}u_{n+1}^{p},\end{displaymath} \ \par

\begin{displaymath} 
Av_{n+1}-\alpha u_{n+1}+\frac{p}{(1-p)\Delta t}v_{n}^{p-1}v_{n+1}=\lambda 
\gamma ^{1-p}p(p+1)v_{n+1}^{p}\end{displaymath} \ \par
and 
\begin{equation} 
\gamma =\displaystyle \left({\displaystyle \frac{\displaystyle 1}{\displaystyle 
(1-p)\Delta tJ_{n}({\hat u},\hat v)}}\right) ^{1/(1-p)}.\label{systeme23}
\end{equation} \ \par
Hence, the numerical scheme admits at least one positive solution.\ \par
\ \par
Before proving the uniqueness of the solution, we first prove the existence 
of bounded supersolutions and  of maximal solutions.\ \par
\ \par
\begin{Lemma} If the hypotheses of the lemma ~(\ref{systeme17}) are satisfied, 
the system ~(\ref{systeme11}), ~(\ref{systeme12}) admits a constant supersolution. 
\ \par
\end{Lemma}
\textsl{Proof: } Let $(C_{n}^{1},C_{n}^{2})\in {\mathbb{R}}_{+}^{2}$; 
$(C_{n}^{1},\ C_{n}^{2})$ will be a supersolution of the system ~(\ref{systeme11}), 
~(\ref{systeme12}) if these constants satisfy the inequalities:\ \par

\begin{equation} 
-\frac{m}{(1-p)\Delta t}(C_{n}^{1})^{p}u_{n}^{m-p}+\frac{m}{(1-p)\Delta 
t}C_{n}^{1}u_{n}^{m-1}\geq \alpha C_{n}^{2},\label{systeme18}
\end{equation} \ \par

\begin{equation} 
-\frac{p}{(1-p)\Delta t}(C_{n}^{2})^{p}+\frac{p}{(1-p)\Delta t}C_{n}^{2}v_{n}^{p-1}\geq 
\alpha C_{n}^{1}.\label{systeme19}
\end{equation} \ \par
We note $x=\frac{C_{n}^{2}}{C_{n}^{1}}$\ \par
The first inequality may be written:\ \par

\begin{displaymath} 
\frac{m}{(1-p)\Delta t}u_{n}^{m-1}-\alpha x\geq \frac{m}{(1-p)\Delta t}(C_{n}^{1})^{p-1}u_{n}^{m-p},\end{displaymath} 
\ \par
hence this inequality may be satisfied only for $x<\frac{m}{(1-p)\alpha \Delta t}\left\Vert{u_{n}}\right\Vert _{\infty 
}^{m-1}.$\ \par
The second inequality may be written:\ \par

\begin{displaymath} 
\frac{p}{(1-p)\Delta t}v_{n}^{p-1}-\frac{\alpha }{x}\geq \frac{p}{(1-p)\Delta 
t}(C_{n}^{2})^{p-1}\end{displaymath} \ \par
and may be satisfied only if  $\displaystyle \frac{\displaystyle 1}{\displaystyle x}<\frac{\displaystyle p}{\displaystyle 
(1-p)\alpha \Delta t}\displaystyle \left\Vert{\displaystyle v_{n}}\right\Vert 
_{\infty }^{p-1}$.\ \par
\ \par
So a necessary condition to obtain inequalities ~(\ref{systeme18}), ~(\ref{systeme19}) 
is that :\ \par
\ \par
$\displaystyle \frac{\displaystyle (1-p)\alpha \Delta t}{\displaystyle p}\displaystyle 
\left\Vert{\displaystyle v_{n}}\right\Vert _{\infty }^{1-p}<\frac{\displaystyle 
m}{\displaystyle (1-p)\alpha \Delta t}\displaystyle \left\Vert{\displaystyle 
u_{n}}\right\Vert _{\infty }^{m-1}$, that is the condition ~(\ref{systeme13}) 
\ \par
\ \par
and $\displaystyle \frac{\displaystyle (1-p)\alpha \Delta t}{\displaystyle p}\displaystyle 
\left\Vert{\displaystyle v_{n}}\right\Vert _{\infty }^{1-p}<x<\frac{\displaystyle 
m}{\displaystyle (1-p)\alpha \Delta t}\displaystyle \left\Vert{\displaystyle 
u_{n}}\right\Vert _{\infty }^{m-1}.$\ \par
\ \par
\ \par
If we choose $\displaystyle C_{n}^{1}=\frac{\displaystyle \displaystyle \left\Vert{\displaystyle u_{n}}\right\Vert 
_{\infty }}{\displaystyle \displaystyle \left({\displaystyle 1-\frac{\displaystyle 
\alpha (1-p)}{\displaystyle m}\Delta tx\displaystyle \left\Vert{\displaystyle 
u_{n}}\right\Vert _{\infty }^{1-m}}\right) ^{1/(1-p)}}$, then we get 
:\ \par
$\displaystyle C_{n}^{2}=\frac{\displaystyle x\displaystyle \left\Vert{\displaystyle 
u_{n}}\right\Vert _{\infty }}{\displaystyle \displaystyle \left({\displaystyle 
1-\frac{\displaystyle \alpha (1-p)}{\displaystyle m}\Delta tx\displaystyle 
\left\Vert{\displaystyle u_{n}}\right\Vert _{\infty }^{1-m}}\right) ^{1/(1-p)}}$.\ 
\par
It remains to prove that $x$ may be chosen in the interval $\displaystyle 
[\frac{\displaystyle (1-p)\alpha \Delta t}{\displaystyle p}\displaystyle 
\left\Vert{\displaystyle v_{n}}\right\Vert _{\infty }^{1-p},\frac{\displaystyle 
m}{\displaystyle (1-p)\alpha \Delta t}\displaystyle \left\Vert{\displaystyle 
u_{n}}\right\Vert _{\infty }^{1-m}]$ such that $C_{n}^{2}$ satisfies 
: $\displaystyle (C_{n}^{2})^{p-1}\leq v_{n}^{p-1}-\frac{\displaystyle \alpha }{\displaystyle 
x}\frac{\displaystyle 1-p}{\displaystyle p}\Delta t$.\ \par
In order to obtain this inequality, the parameter $x$ must satisfy:\ 
\par

\begin{equation} 
\frac{\left\Vert{v_{n}}\right\Vert _{\infty }}{\left({1-\alpha \frac{1-p}{p}\frac{\Delta 
t}{x}\left\Vert{v_{n}}\right\Vert _{\infty }^{1-p}}\right) ^{1/(1-p)}}\leq 
\frac{x\left\Vert{u_{n}}\right\Vert _{\infty }}{\left({1-\alpha \frac{1-p}{m}\Delta 
tx\left\Vert{u_{n}}\right\Vert _{\infty }^{1-m}}\right) ^{1/(1-p)}}\label{systeme20}
\end{equation} \ \par
Let us denote $a=\frac{1-p}{p}\Delta t\alpha \left\Vert{v_{n}}\right\Vert _{\infty }^{1-p},$ 
$b=\frac{m}{(1-p)\alpha \Delta t}\left\Vert{u_{n}}\right\Vert _{\infty 
}^{m-1}$;\ \par
\ \par
we have : $a<x<b$ and the condition ~(\ref{systeme20}) may be written: 
\ \par
{\hskip 1.2em}$\displaystyle \frac{\displaystyle p}{\displaystyle (1-p)\alpha \Delta t}ab^{(1-p)/(1-m)}\displaystyle 
\left({\displaystyle 1-\frac{\displaystyle x}{\displaystyle b}}\right) 
\leq \displaystyle \left({\displaystyle \frac{\displaystyle m}{\displaystyle 
(1-p)\alpha \Delta t}}\right) ^{(1-p)/(1-m)}\displaystyle \left({\displaystyle 
1-\frac{\displaystyle a}{\displaystyle x}}\right) x^{1-p}$.\ \par
\ \par
If we define the function $f$ by\ \par

\begin{displaymath} 
f(x)=\frac{p}{(1-p)\alpha \Delta t}ab^{(m-p)/(1-m)}x^{p}(x-b)+\left({\frac{m}{(1-p)\alpha 
\Delta t}}\right) ^{(1-p)/(1-m)}(x-a)\end{displaymath} \ \par
the condition ~(\ref{systeme20}) becomes: $f(x)\geq 0$.\ \par
\ \par
The function $f$ satisfies $f(a)<0$ and $f(b)>0$; so there exists $x_{0}\in ]a,b[$ 
such that $f(x_{0})=0$\ \par
and the couple \ \par
$\displaystyle \displaystyle \left({\displaystyle C_{n}^{1}=\frac{\displaystyle \displaystyle 
\left\Vert{\displaystyle u_{n}}\right\Vert _{\infty }}{\displaystyle 
\displaystyle \left({\displaystyle 1-\alpha \frac{\displaystyle (1-p)}{\displaystyle 
m}\Delta tx_{0}\displaystyle \left\Vert{\displaystyle u_{n}}\right\Vert 
_{\infty }^{1-m}}\right) ^{1/(1-p)}},C_{n}^{2}=\frac{\displaystyle \displaystyle 
\left\Vert{\displaystyle v_{n}}\right\Vert _{\infty }}{\displaystyle 
\displaystyle \left({\displaystyle 1-\alpha \frac{\displaystyle 1-p}{\displaystyle 
p}\frac{\displaystyle \Delta t}{\displaystyle x_{0}}\displaystyle \left\Vert{\displaystyle 
v_{n}}\right\Vert _{\infty }^{1-p}}\right) ^{1/(1-p)}}}\right) $\ \par
\ \par
is a supersolution of the system ~(\ref{systeme11}), ~(\ref{systeme12}). 
\ \par
\ \par
\begin{Lemma} System ~(\ref{systeme11}),~(\ref{systeme12}) has a maximal 
solution $(\overline{u},\ \overline{v})$ and any solution $(u,v)$ satisfies: 
$0\leq u\leq \overline{u},\ 0\leq v\leq \overline{v}$.\ \par
\end{Lemma}
\textsl{Proof: }We use the same method as Keller in ~\cite{Keller}.\textsl{ 
}We consider the sequences defined by:$\ u_{n+1,0}=C_{n}^{1},\ v_{n+1,0}=C_{n}^{2},$\ 
\par

\begin{displaymath} 
Au_{n+1,j+1}+\frac{m}{(1-p)\Delta t}u_{n}^{m-1}u_{n+1,j+1}=\alpha v_{n+1,j}+\frac{m}{(1-p)\Delta 
t}u_{n+1,j}^{p}u_{n}^{m-p},\end{displaymath} \ \par

\begin{displaymath} 
Av_{n+1,j+1}+\frac{p}{(1-p)\Delta t}v_{n}^{p-1}v_{n+1,j+1}=\alpha u_{n+1,j}+\frac{p}{(1-p)\Delta 
t}v_{n+1,j}^{p}.\end{displaymath} \ \par
We get:
\begin{displaymath} 
A(u_{n+1,1}-u_{n+1,0})+\frac{m}{(1-p)\Delta t}u_{n}^{m-1}(u_{n+1,1}-u_{n+1,0})=\alpha 
C_{n}^{2}+\frac{m}{(1-p)\Delta t}u_{n}^{m-p}\left({(C_{n}^{1})^{p}-u_{n}^{p-1}C_{n}^{1}}\right) 
.\end{displaymath} \ \par
The second member of this equality is negative; we deduce from the maximum 
principle that: $u_{n+1,1}\leq u_{n+1,0}.$ In the same manner, we get: 
$v_{n+1,1}\leq v_{n+1,0}$. We prove recurently that the sequences $(u_{n+1,j})_{j\geq 0}$ 
and $(v_{n+1,j)})_{j\geq 0}$ are decreasing; in fact, we have:\ \par
\ \par
$\displaystyle A(u_{n+1,j+1}-u_{n+1,j})+\frac{\displaystyle m}{\displaystyle (1-p)\Delta 
t}u_{n}^{m-1}(u_{n+1,j+1}-u_{n+1,j})=$\ \par
$\displaystyle \alpha (v_{n+1,j}-v_{n+1,j-1})+\frac{\displaystyle m}{\displaystyle (1-p)\Delta 
t}u_{n}^{m-p}(u_{n+1,j}^{p}-u_{n+1,j-1}^{p})$\ \par
and the second member is negative from the recurrence hypothesis.\ \par
\ \par
We deduce: $u_{n+1,j+1}\leq u_{n+1,j}$. Similarly, we get: $v_{n+1,j+1}\leq v_{n+1,j}.$\ 
\par
\ \par
Since the two sequences $(u_{n+1,j})_{j\geq 0}$ and $(v_{n+1,j})_{j\geq 0}$ 
are nonnegative, they converge to $\overline{u}$ and  $\overline{v}$ 
 and taking the limit when $j{\longrightarrow}+\infty $, we obtain:\ 
\par

\begin{displaymath} 
A\overline{u}+\frac{m}{(1-p)\Delta t}u_{n}^{m-1}\overline{u}=\alpha \overline{v}+\frac{m}{(1-p)\Delta 
t}u_{n}^{m-p}\overline{u}^{p}\end{displaymath} \ \par

\begin{displaymath} 
A\overline{v}+\frac{p}{(1-p)\Delta t}v_{n}^{p-1}\overline{v}=\alpha \overline{u}+\frac{p}{(1-p)\Delta 
t}\overline{v}^{p}.\end{displaymath} \ \par
So, the functions $\overline{u}$ and $\overline{v}$ are solutions of 
the system ~(\ref{systeme11}), ~(\ref{systeme12}). \ \par
It remains to prove that any solution $(u,v)$ satisfies: $0\leq u\leq \overline{u},\ 0\leq v\leq \overline{v}$.\ 
\par
Let $(u,v)$ a solution of system ~(\ref{systeme11}), ~(\ref{systeme12}), 
we have: $0\leq u\leq C_{n}^{1},\ 0\leq v\leq C_{n}^{2}$; we have the 
equalities:\ \par

\begin{displaymath} 
A(u-u_{n+1,j+1})+\frac{m}{(1-p)\Delta t}u_{n}^{m-1}(u-u_{n+1,j+1})=\alpha 
(v-v_{n+1,j})+\frac{m}{(1-p)\Delta t}u_{n}^{m-p}(u^{p}-u_{n+1,j}^{p})\end{displaymath} 
\ \par

\begin{displaymath} 
A(v-v_{n+1,j+1})+\frac{p}{(1-p)\Delta t}v_{n}^{p-1}(v-v_{n+1,j+1})=\alpha 
(u-u_{n+1,j})+\frac{p}{(1-p)\Delta t}(v^{p}-v_{n+1,j}^{p}).\end{displaymath} 
\ \par
For $j=0,$ the second member of these inequalities is negative; then 
we get $u\leq u_{n+1,1}$ , $v\leq v_{n+1,1}$ and recurrently, we obtain 
$u\leq u_{n+1,j},\ v\leq v_{n+1,j}$ for any $j\geq 0$. It results: $u\leq \overline{u},\ v\leq \overline{v}$.\ 
\par
\ \par
\begin{Theorem} If the functions $u_{n}$ and $v_{n}$ are positive, continuous 
in $\overline{\Omega }$ and satisfy the condition ~(\ref{systeme13}), 
then system ~(\ref{systeme11}), ~(\ref{systeme13}) has a unique positive 
solution.\ \par
\end{Theorem}
\textsl{Proof: } From the previous lemmas, we know that the system admits 
at least one positive solution and that any solution $(u,v)$ satisfies 
$0\leq u\leq \overline{u}$,  $0\leq v\leq \overline{v}$.\ \par
We get: 
\begin{displaymath} 
\int_{\Omega }^{}{(Au\ \overline{u}-A\overline{u}\ u)dx}=0=\frac{m}{(1-p)\Delta 
t}\int_{\Omega }^{}{u_{n}^{m-p}u\ \overline{u}(u^{p-1}\ \ -\overline{u}^{p-1})dx}+\alpha 
\int_{\Omega }^{}{(v\overline{u}-\overline{v}u)dx}.\end{displaymath} 
\ \par
Similarly, we have: \ \par

\begin{displaymath} 
\ \frac{p}{(1-p)\Delta t}\int_{\Omega }^{}{v\ \overline{v}\ (v^{p-1}-\overline{v}^{p-1})dx}+\alpha 
\int_{\Omega }^{}{(u\overline{v}-\overline{u}v)dx}=0\end{displaymath} 
\ \par
We deduce from these equalities that $\int_{\Omega }^{}{(u\overline{v}\ -\ v\overline{u})dx}=0$ 
and then $u=\overline{u},\ v=\overline{v}$.\ \par
\ \par
\begin{Theorem} The numerical solution exists at least until the time 
\ \par

\begin{displaymath} 
T_{1}=\min \left({\frac{m}{\alpha (1-p)}\lambda _{0}^{m-1},\frac{p}{\alpha 
(1-p)}\lambda _{0}^{p-1}}\right) \end{displaymath} \ \par
with $\lambda _{0}=\max \left({\left\Vert{u_{0}}\right\Vert _{\infty },\ \left\Vert{v_{0}}\right\Vert 
_{\infty }}\right) $.\label{systeme34}\ \par
\end{Theorem}
\ \par
\textsl{Proof:} We prove recurently that the solution $(u_{n},v_{n)}$ 
satisfy the inequality: $\left\Vert{u_{n}}\right\Vert _{\infty },\ \left\Vert{v_{n}}\right\Vert 
_{\infty }\leq \phi _{n}$ \ \par
where $\phi _{n}$ is defined by  $\displaystyle \phi _{0}=\max \displaystyle \left({\displaystyle \frac{\displaystyle 
\alpha (1-p)}{\displaystyle m}\lambda _{0}^{m-1},\ \frac{\displaystyle 
\alpha (1-p)}{\displaystyle p}\lambda _{0}^{p-1}}\right) $ , $\displaystyle 
\phi _{n}=\frac{\displaystyle \lambda _{0}}{\displaystyle (1-t_{n}\phi 
_{0})^{1/(1-p)}}.$\ \par
If this inequality is satisfied at the time level $t_{n}=n\Delta t,\ $if 
$t_{n+1}\phi _{0}\leq 1$, the inequality ~(\ref{systeme13}) is verified 
and the solution exists at the time level $t_{n+1}.$\ \par
The quantity $\phi _{n+1}$ will be a supersolution of the system ~(\ref{systeme11}), 
~(\ref{systeme12}), if we have the two inequalities:\ \par
\ \par

\begin{displaymath} 
-\frac{\displaystyle m}{\displaystyle (1-p)\Delta t}\phi _{n+1}^{p}u_{n}^{m-p}+\frac{\displaystyle 
m}{\displaystyle (1-p)\Delta t}\phi _{n+1}u_{n}^{m-1}\geq \alpha \phi 
_{n+1},\end{displaymath} \ \par

\begin{displaymath} 
-\frac{\displaystyle p}{\displaystyle (1-p)\Delta t}\phi _{n+1}^{p}+\frac{\displaystyle 
p}{\displaystyle (1-p)\Delta t}\phi _{n+1}v_{n}^{p-1}\geq \alpha \phi 
_{n+1}.\end{displaymath} \ \par
This may be written:\ \par
$\displaystyle \phi _{n+1}^{1-p}\geq \max \displaystyle \left({\displaystyle \frac{\displaystyle 
\displaystyle \left\Vert{\displaystyle u_{n}}\right\Vert _{\infty }^{1-p}}{\displaystyle 
1-\alpha \frac{\displaystyle (1-p)}{\displaystyle m}\Delta t\displaystyle 
\left\Vert{\displaystyle u_{n}}\right\Vert _{\infty }^{1-m}},\ \frac{\displaystyle 
\displaystyle \left\Vert{\displaystyle v_{n}}\right\Vert _{\infty }^{1-p}}{\displaystyle 
1-\alpha \frac{\displaystyle (1-p)}{\displaystyle p}\Delta t\displaystyle 
\left\Vert{\displaystyle v_{n}}\right\Vert _{\infty }^{1-p}}}\right) 
$.\ \par
But, from the recurrence hypothese, we get \ \par

\begin{displaymath} 
\frac{\displaystyle \displaystyle \left\Vert{\displaystyle u_{n}}\right\Vert 
_{\infty }^{1-p}}{\displaystyle 1-\alpha \frac{\displaystyle (1-p)}{\displaystyle 
m}\Delta t\displaystyle \left\Vert{\displaystyle u_{n}}\right\Vert _{\infty 
}^{1-m}}\leq \frac{\displaystyle \lambda _{0}^{1-p}}{\displaystyle (1-t_{n}\phi 
_{0})\displaystyle \left({\displaystyle 1-\alpha \frac{\displaystyle 
1-p}{\displaystyle m}\Delta t\frac{\displaystyle \lambda _{0}^{1-m}}{\displaystyle 
(1-t_{n}\phi _{0})^{(1-m)/(1-p)}}}\right) }\end{displaymath} \ \par
and it is easy to see that this quantitiy is bounded by $\phi _{n+1}.$\ 
\par
In an analogous manner, we obtain that $\displaystyle \frac{\displaystyle \displaystyle \left\Vert{\displaystyle v_{n}}\right\Vert 
_{\infty }^{1-p}}{\displaystyle 1-\alpha \frac{\displaystyle 1-p}{\displaystyle 
p}\Delta t\displaystyle \left\Vert{\displaystyle v_{n}}\right\Vert _{\infty 
}^{1-p}}$ is bounded by $\phi _{n+1}.$\ \par
So, the solution at the time level $t_{n+1}\ $satisfies: $\left\Vert{u_{n+1}}\right\Vert _{\infty },\ \left\Vert{v_{n+1}}\right\Vert 
_{\infty }\leq \phi _{n+1}$ and the numerical solution exists exists 
during a positive time interval.\ \par
\ \par
\ \par
\section{Properties of the numerical scheme}
\setcounter{equation}{0}In this section, we prove that if $\alpha >\lambda _{1}$, 
the numerical solution blows up in finite time.\ \par
We define the functional $\psi _{n}$ and $F_{n}$ by: 
\begin{displaymath} 
\psi _{n}(u,v)=\left({\int_{\Omega }^{}{(mu_{n}^{m-p}u^{p+1}+pv^{p+1})dx}}\right) 
^{1/(p+1)}\end{displaymath} \ \par
and 
\begin{displaymath} 
F_{n}(u,v)=\frac{J(u,v)}{(\psi _{n}(u,v))^{2}}.\end{displaymath} \ \par
\begin{Lemma} The sequence $(F_{n}(u_{n},v_{n}))_{n\geq 0}$ is nonincreasing.\ 
\par
\end{Lemma}
\textsl{Proof:  }Since $u_{n+1}=\gamma {\hat u}$ and $v_{n+1}=\gamma \hat v,\ ({\hat u},\hat v)\in K$, 
we get $\psi _{n}(u_{n+1},v_{n+1})=\gamma $  and $F_{n}(u_{n+1},v_{n+1})=J({\hat u},\hat v)$.\ 
\par
Besides from ~(\ref{systeme24}), we have  $\displaystyle J({\hat u},\hat v)=J_{n}({\hat u},\hat v)-\frac{\displaystyle 1}{\displaystyle 
(1-p)\Delta t}\displaystyle \int_{\Omega }^{}{(mu_{n}^{m-1}{\hat u}^{2}+pv_{n}^{p-1}\hat 
v^{2})dx}.$\ \par
\ \par
Hence, we get : $\displaystyle J({\hat u},\hat v)\leq \frac{\displaystyle J_{n}(u_{n},v_{n})}{\displaystyle 
\psi _{n}^{2}(u_{n},v_{n})}-\frac{\displaystyle 1}{\displaystyle (1-p)\Delta 
t}\frac{\displaystyle \displaystyle \int_{\Omega }^{}{(mu_{n}^{m-1}u_{n+1}^{2}+pv_{n}^{p-1}v_{n+1}^{2})dx}}{\displaystyle 
\psi _{n}^{2}(u_{n+1},v_{n+1})}.$\ \par
\ \par
In addition, we have the equality:\ \par
$\displaystyle J_{n}(u_{n},v_{n})=J(u_{n},v_{n})+\frac{\displaystyle 1}{\displaystyle 
(1-p)\Delta t}(\psi _{n}(u_{n},v_{n}))^{p+1}$.\ \par
We deduce:\ \par
$\displaystyle J({\hat u},\hat v)\leq \frac{\displaystyle J(u_{n},v_{n})}{\displaystyle 
\psi _{n}^{2}(u_{n},v_{n})}+\frac{\displaystyle 1}{\displaystyle (1-p)\Delta 
t}\displaystyle \left({\displaystyle (\psi _{n}(u_{n},v_{n}))^{p-1}-\frac{\displaystyle 
\displaystyle \int_{\Omega }^{}{(mu_{n}^{m-1}u_{n+1}^{2}+pv_{n}^{p-1}v_{n+1}^{2})dx}}{\displaystyle 
\psi _{n}^{2}(u_{n+1},v_{n+1})}}\right) .$\ \par
\ \par
By the H{\"o}lder inequality, we have at once:\ \par
$\displaystyle \displaystyle \int_{\Omega }^{}{u_{n}^{m-p}u_{n+1}^{p+1}dx}\leq \displaystyle 
\left({\displaystyle \displaystyle \int_{\Omega }^{}{u_{n+1}^{2}u_{n}^{m-1}dx}}\right) 
^{(p+1)/2}\displaystyle \left({\displaystyle \displaystyle \int_{\Omega 
}^{}{u_{n}^{m+1}dx}}\right) ^{(1-p)/2},$\ \par
\ \par
$\displaystyle \displaystyle \int_{\Omega }^{}{v_{n+1}^{p+1}dx}\leq \displaystyle \left({\displaystyle 
\displaystyle \int_{\Omega }^{}{v_{n+1}^{2}v_{n}^{p-1}dx}}\right) ^{(p+1)/2}\displaystyle 
\left({\displaystyle \displaystyle \int_{\Omega }^{}{v_{n}^{p+1}dx}}\right) 
^{(1-p)/2}.$ \ \par
\ \par
Hence, we get:\ \par
\ \par
$\displaystyle \psi _{n}^{p+1}(u_{n+1},v_{n+1})\leq m\displaystyle \left({\displaystyle 
\displaystyle \int_{\Omega }^{}{u_{n+1}^{2}u_{n}^{m-1}dx}}\right) ^{(p+1)/2}\displaystyle 
\left({\displaystyle \displaystyle \int_{\Omega }^{}{u_{n}^{m+1}dx}}\right) 
^{(1-p)/2}$\ \par
$\displaystyle \ +p\displaystyle \left({\displaystyle \displaystyle \int_{\Omega }^{}{v_{n+1}^{2}v_{n}^{p-1}dx}}\right) 
^{(p+1)/2}\displaystyle \left({\displaystyle \displaystyle \int_{\Omega 
}^{}{v_{n}^{p+1}dx}}\right) ^{(1-p)/2}\ $\ \par
and \ \par

\begin{equation} 
\psi _{n}^{2}(u_{n+1},v_{n+1})\leq \ \left({m\int_{\Omega }^{}{u_{n+1}^{2}u_{n}^{m-1}dx}+p\int_{\Omega 
}^{}{v_{n+1}^{2}v_{n}^{p-1}dx}}\right) \psi _{n}^{1-p}(u_{n},v_{n}).\label{systeme25}
\end{equation} \ \par
We deduce: $\displaystyle J({\hat u},\hat v)\leq \frac{\displaystyle J(u_{n},v_{n})}{\displaystyle 
\psi _{n}^{2}(u_{n},v_{n})},\ $that is $F_{n}(u_{n+1},v_{n+1})\leq F_{n}(u_{n},v_{n}).$\ 
\par
\ \par
\begin{Lemma} For $n\geq 0$, we have the estimate:\ \par

\begin{equation} 
(1-p)\Delta tF_{n}(u_{n+1},v_{n+1})\leq \psi _{n}^{p-1}(u_{n+1},v_{n+1})-\psi 
_{n}^{p-1}(u_{n},v_{n})\leq (1-p)\Delta tF_{n}(u_{n},v_{n}).\label{systeme26}
\end{equation} \ \par
\end{Lemma}
\textsl{Proof: }$\iota $\textsl{) }We prove first the right inequality. 
We have:  $\psi (u_{n+1},v_{n+1})=\gamma $;  from ~(\ref{systeme23}), 
we obtain\ \par

\begin{displaymath} 
\psi _{n}(u_{n+1},v_{n+1})=\left({(1-p)\Delta tJ_{n}({\hat u},\hat v)}\right) 
^{1/(p-1)}\end{displaymath} \ \par
and 
\begin{displaymath} 
\psi _{n}^{p-1}(u_{n+1},v_{n+1})\leq (1-p)\Delta t\frac{J_{n}(u_{n},v_{n})}{\psi 
_{n}^{2}(u_{n},v_{n})},\end{displaymath}  that is\ \par
$\psi _{n}^{p-1}(u_{n+1},v_{n+1})\leq (1-p)\Delta t\left({F_{n}(u_{n},v_{n})+\frac{1}{(1-p)\Delta 
t}\psi _{n}^{p-1}(u_{n},v_{n})}\right) .$\ \par
\ \par
$\iota \iota $)Multiplying ~(\ref{systeme11}) by $u_{n+1}$ and ~(\ref{systeme12}) 
by $v_{n+1}$ and integrating on $\Omega $, we get:\ \par

\begin{displaymath} 
\frac{\displaystyle m}{\displaystyle (1-p)\Delta t}\displaystyle \int_{\Omega 
}^{}{u_{n+1}^{2}u_{n}^{m-1}dx}-\frac{\displaystyle m}{\displaystyle (1-p)\Delta 
t}\displaystyle \int_{\Omega }^{}{u_{n+1}^{p+1}u_{n}^{m-p}dx}+\displaystyle 
\int_{\Omega }^{}{\displaystyle \left({\displaystyle \displaystyle \left\vert{\displaystyle 
\nabla u_{n+1}}\right\vert ^{2}-\alpha u_{n+1}v_{n+1}}\right) dx}=0,\end{displaymath} 
\ \par

\begin{displaymath} 
\frac{\displaystyle p}{\displaystyle (1-p)\Delta t}\displaystyle \int_{\Omega 
}^{}{v_{n+1}^{2}v_{n}^{p-1}dx}-\frac{\displaystyle p}{\displaystyle (1-p)\Delta 
t}\displaystyle \int_{\Omega }^{}{v_{n+1}^{p+1}dx}+\displaystyle \int_{\Omega 
}^{}{\displaystyle \left({\displaystyle \displaystyle \left\vert{\displaystyle 
\nabla v_{n+1}}\right\vert ^{2}-\alpha u_{n+1}v_{n+1}}\right) dx}=0.\end{displaymath} 
\ \par
Hence , we get :\ \par

\begin{displaymath} 
\frac{1}{(1-p)\Delta t}\int_{\Omega }^{}{(mu_{n+1}^{2}u_{n}^{m-1}+pv_{n+1}^{2}v_{n}^{p-1})dx}-\frac{1}{(1-p)\Delta 
t}\psi _{n}^{p+1}(u_{n+1},v_{n+1})+J(u_{n+1},v_{n+1})=0.\end{displaymath} 
\ \par
By using ~(\ref{systeme25}), we deduce:\ \par

\begin{displaymath} 
\psi _{n}^{p-1}(u_{n},v_{n})-\psi _{n}^{p-1}(u_{n+1},v_{n+1})+F(u_{n+1},v_{n+1})\leq 
0.\end{displaymath} \ \par
This concludes the proof.\ \par
\ \par
\begin{Lemma} The sequence $(J(u_{n},v_{n}))_{n\geq 0}$ is nonincreasing\ 
\par
\end{Lemma}
\textsl{Proof:} In ~\cite{mnLeRoux1} , we have proved the inequality:\ 
\par

\begin{displaymath} 
\forall a,b\in {\mathbb{R}}^{+},a^{p-1}(b-a)^{2}\leq a^{p+1}-b^{p+1}-\frac{1+p}{1-p}b^{2}(b^{p-1}-a^{p-1)}.\end{displaymath} 
\ \par
We deduce from this inequality:\ \par

\begin{displaymath} 
\int_{\Omega }^{}{(u_{n+1}-u_{n})^{2}u_{n}^{m-1}dx}\leq \int_{\Omega }^{}{u_{n}^{m+1}dx}-\int_{\Omega 
}^{}{u_{n+1}^{p+1}u_{n}^{m-p}dx}-\frac{1+p}{1-p}\int_{\Omega }^{}{u_{n+1}^{2}(u_{n+1}^{p-1}-u_{n}^{p-1})u_{n}^{m-p}dx}\end{displaymath} 
\ \par
and\ \par

\begin{displaymath} 
\int_{\Omega }^{}{(v_{n+1}-v_{n})^{2}v_{n}^{p-1}dx}\leq \int_{\Omega }^{}{v_{n}^{p+1}dx}-\int_{\Omega 
}^{}{v_{n+1}^{p+1}dx}-\frac{1+p}{1-p}\int_{\Omega }^{}{v_{n+1}^{2}(v_{n+1}^{p-1}-v_{n}^{p-1})dx}.\end{displaymath} 
\ \par
Since $u_{n+1}$ and $v_{n+1}$ are solutions of ~(\ref{systeme11}), ~(\ref{systeme12}), 
we get:\ \par

\begin{displaymath} 
\int_{\Omega }^{}{u_{n+1}^{2}(u_{n+1}^{p-1}-u_{n}^{p-1})u_{n}^{m-p}dx}=\frac{(1-p)\Delta 
t}{m}\int_{\Omega }^{}{(\left\vert{\nabla u_{n+1}}\right\vert ^{2}-\alpha 
u_{n+1}v_{n+1})dx}\end{displaymath} \ \par
and\ \par

\begin{displaymath} 
\int_{\Omega }^{}{v_{n+1}^{2}(v_{n+1}^{p-1}-v_{n}^{p-1})dx}=\frac{(1-p)\Delta 
t}{p}\int_{\Omega }^{}{(\left\vert{\nabla v_{n+1}}\right\vert ^{2}-\alpha 
u_{n+1}v_{n+1})dx}.\end{displaymath} \ \par
So, we obtain:
\begin{displaymath} 
m\int_{\Omega }^{}{(u_{n+1}-u_{n})^{2}u_{n}^{m-1}dx}+p\int_{\Omega }^{}{(v_{n+1}-v_{n})^{2}v_{n}^{p-1}dx}\end{displaymath} 
\ \par

\begin{displaymath} 
\leq \psi _{n}^{p+1}(u_{n},v_{n})-\psi _{n}^{p+1}(u_{n+1},v_{n+1})-(1+p)\Delta 
tJ(u_{n+1},v_{n+1}).\end{displaymath} \ \par
From  the inequality ~\cite{mnLeRoux1} :$\forall a,b\in {\mathbb{R}}^{+},\ a^{p+1}-b^{p+1}\leq \frac{1+p}{1-p}\ 
a^{2}(b^{p-1}-a^{p-1}),$\ \par
we deduce: 
\begin{displaymath} 
\psi _{n}^{p+1}(u_{n},v_{n})-\psi _{n}^{p+1}(u_{n+1},v_{n+1})\leq \frac{1+p}{1-p}\psi 
_{n}^{2}(u_{n},v_{n})(\psi _{n}^{p-1}(u_{n+1},v_{n+1})-\psi _{n}^{p-1}(u_{n},v_{n}))\end{displaymath} 
and by using ~(\ref{systeme26}), we get:\ \par

\begin{equation} 
\psi _{n}^{p+1}(u_{n},v_{n})-\psi _{n}^{p+1}(u_{n+1},v_{n+1})\leq (1+p)\Delta 
tJ(u_{n},v_{n}).\label{systeme27}
\end{equation} \ \par
So, we get \ \par

\begin{displaymath} 
m\int_{\Omega }^{}{(u_{n+1}-u_{n})^{2}u_{n}^{m-1}dx}+p\int_{\Omega }^{}{(v_{n+1}-v_{n})^{2}v_{n}^{p-1}dx}\leq 
(1+p)\Delta t\left({J(u_{n},v_{n})-J(u_{n+1},v_{n+1})}\right) .\end{displaymath} 
We deduce: $J(u_{n},v_{n})\geq J(u_{n+1},v_{n+1})$.\ \par
\ \par
We note 
\begin{displaymath} 
\Phi _{n}=\displaystyle \int_{\Omega }^{}{\displaystyle \left({\displaystyle 
\frac{\displaystyle m}{\displaystyle m+1}u_{n}^{m+1}+\frac{\displaystyle 
p}{\displaystyle p+1}v_{n}^{p+1}}\right) dx}\end{displaymath} \ \par
\ \par
\begin{Lemma} If $J(u_{0},v_{0})<0,$ the sequence $(\Phi _{n})_{n\geq 0}$ 
is increasing.\ \par
\end{Lemma}
\textsl{Proof:} We have the equality :\ \par

\begin{equation} 
\psi _{n}^{p+1}(u_{n},v_{n})=\displaystyle \int_{\Omega }^{}{(mu_{n}^{m+1}+pv_{n}^{p+1})dx}=(p+1)\Phi 
_{n}+\frac{\displaystyle m(m-p)}{\displaystyle m+1}\displaystyle \int_{\Omega 
}^{}{u_{n}^{m+1}dx}.\label{systeme28}
\end{equation} \ \par
Besides, we get:\ \par

\begin{equation} 
\psi _{n}^{p+1}(u_{n+1},v_{n+1})=\displaystyle \int_{\Omega }^{}{(mu_{n}^{m-p}u_{n+1}^{p+1}+pv_{n+1}^{p+1})dx}\leq 
(p+1)\Phi _{n+1}+\frac{\displaystyle m(m-p)}{\displaystyle m+1}\displaystyle 
\int_{\Omega }^{}{u_{n}^{m+1}dx}.\label{systeme29}
\end{equation} \ \par
and we deduce:
\begin{displaymath} 
\Phi _{n}-\Phi _{n+1}\leq \frac{1}{p+1}\left({\psi _{n}^{p+1}(u_{n},v_{n})-\psi 
_{n}^{p+1}(u_{n+1},v_{n+1})}\right) .\end{displaymath} \ \par
By using ~(\ref{systeme27}) , we obtain $\Phi _{n}-\Phi _{n+1}\leq \Delta tJ(u_{n},v_{n}).$\ 
\par
If $J(u_{0},v_{0})<0,$ since the sequence $(J(u_{n},v_{n}))_{n\geq 0}$ 
is nonincreasing, we deduce that the sequence $(\Phi _{n})_{n\geq 0}$ 
is increasing.\ \par
\ \par
\begin{Lemma} If $J(u_{0},v_{0})<0$, for $n\geq 0,$ we have the inequality: 
\ \par

\begin{equation} 
\Phi _{n}^{2/(p+1)}\left({\Phi _{n+1}^{(p-1)/(p+1)}-\Phi _{n}^{(p-1)/(p+1)}}\right) 
\leq \frac{1}{m+1}\psi _{n}^{2}(u_{n},v_{n})\left({\psi _{n}^{p-1}(u_{n+1},v_{n+1})-\psi 
_{n}^{p-1}(u_{n},v_{n})}\right) .\label{systeme30}
\end{equation} \ \par
\end{Lemma}
\textsl{Proof: }This inequality may be written:\ \par

\begin{displaymath} 
\psi _{n}^{p+1}(u_{n},v_{n})+(m+1)\Phi _{n}^{2/(p+1)}\Phi _{n+1}^{(p-1)/(p+1)}\leq 
\psi _{n}^{2}(u_{n},v_{n})\psi _{n}^{p-1}(u_{n+1},v_{n+1})+(m+1)\Phi 
_{n}.\end{displaymath} \ \par
By using ~(\ref{systeme28}) and ~(\ref{systeme29}), we obtain that a 
sufficient condition to satify this inequality is:\ \par

\begin{displaymath} 
(p+1)\Phi _{n}+(m-p)\mu _{n}+(m+1)\Phi _{n}^{2/(p+1)}\Phi _{n+1}^{(p-1)/(p+1)}\end{displaymath} 
\ \par

\begin{displaymath} 
\leq \left({(p+1)\Phi _{n}+(m-p)\mu _{n}}\right) ^{2/(p+1)}\left({(p+1)\Phi 
_{n+1}+(m-p)\mu _{n}}\right) ^{(p-1)/(p+1)}+(m+1)\Phi _{n}\end{displaymath} 
\ \par
with $\displaystyle \mu _{n}=\frac{\displaystyle m}{\displaystyle m+1}\displaystyle \int_{\Omega 
}^{}{u_{n}^{m+1}dx}$.\ \par
If $J(u_{0},v_{0})<0,$ the sequence $(\Phi _{n})_{n\geq 0}$ is increasing, 
 so we get: $\displaystyle \frac{\displaystyle (p+1)\Phi _{n}+(m-p)\mu _{n}}{\displaystyle (p+1)\Phi 
_{n+1}+(m-p)\mu _{n}}\geq \frac{\displaystyle \Phi _{n}}{\displaystyle 
\Phi _{n+1}}.$\ \par
Hence, it is sufficient to prove: \ \par

\begin{displaymath} 
(p+1)\Phi _{n}+(m-p)\mu _{n}+(m+1)\Phi _{n}^{2/(p+1)}\Phi _{n+1}^{(p-1)/(p+1)}\end{displaymath} 
\ \par

\begin{displaymath} 
\leq \left({(p+1)\Phi _{n}+(m-p)\mu _{n})}\right) \Phi _{n+1}^{(p-1)/(p+1)}\Phi 
_{n}^{(1-p)/(1+p)}+(m+1)\Phi _{n}\end{displaymath} \ \par
that is $\mu _{n}(\Phi _{n+1}^{(1-p)/(1+p)}-\Phi _{n}^{(1-p)/(1+p)})\leq \Phi _{n}(\Phi 
_{n+1}^{(1-p)/(1+p)}-\Phi _{n}^{(1-p)/(1+p)})$\ \par
\ \par
and this inequality is satisfied since $\mu _{n}\leq \Phi _{n}$.\ \par
\ \par
\begin{Lemma} If $J(u_{0},v_{0})<0$, we have the estimate:\ \par

\begin{equation} 
\Phi _{n}^{2/(p+1)}\left({\Phi _{n+1}^{(p-1)/(p+1)}-\Phi _{n}^{(p-1)/(p+1)}}\right) 
\leq \frac{1-p}{1+m}\Delta tJ(u_{n},v_{n})\label{systeme31}
\end{equation} \ \par
\end{Lemma}
\textsl{Proof: }We deduce the estimate immediately  from ~(\ref{systeme26}) 
and  ~(\ref{systeme30})\ \par
\ \par
\begin{Theorem} If $J(u_{0},v_{0})<0$, the numerical solution blows up 
in a finite time T$_{*}$ such that \ \par

\begin{displaymath} 
T_{*}<\frac{1+m}{1-p}\ \frac{\int_{\Omega }^{}{\left({\frac{m}{m+1}u_{0}^{m+1}+\frac{p}{p+1}v_{0}^{p+1}}\right) 
dx}}{-J(u_{0},v_{0})}.\end{displaymath} \ \par
\ \par
\end{Theorem}
\textsl{Proof: }From ~(\ref{systeme31}), we get $\Phi _{n+1}^{(p-1)/(p+1)}\leq \Phi _{n}^{(p-1)/(p+1)}+\frac{1-p}{1+m}\Delta 
t\Phi _{n}^{-2/(p+1)}J(u_{n},v_{n})$\ \par
and since the sequence $(J(u_{n},v_{n}))_{n\geq 0}\ $is decreasing and 
the sequence $(\Phi _{n})$ increasing, we get :
\begin{displaymath} 
\Phi _{n+1}^{(p-1)/(p+1)}\leq \Phi _{0}^{(p-1)/(p+1)}\left({1+\frac{1-p}{1+m}t_{n}\Phi 
_{0}^{-1}J(u_{0},v_{0})}\right) \end{displaymath} \ \par
and we deduce the estimate.\ \par
\ \par
\begin{Remark} In the case $p=m,$ we obtain the same bound for the numerical 
blow-up time and for the blow-up time of the exact solution. In the case 
$p<m,$ the bound obtained for the numerical blow-up time is inferior 
to the estimate obtained for the blow-up time of the exact solution \ 
\par
\end{Remark}
\ \par
\ \par
\section{The case $p=m$}
\setcounter{equation}{0}\ \par
In the case $p=m$, the functionals $\psi _{n}$ and $F_{n}$ are independent 
of $n$. We shall note them respectively $\psi $ and $F$:\ \par
$\psi (u,v)=\left({\int_{\Omega }^{}{m(u^{m+1}+v^{m+1})dx}}\right) ^{\frac{1}{m+1}}$ 
and $\displaystyle F(u,v)=\frac{\displaystyle J(u,v)}{\displaystyle \psi ^{2}(u,v)}$.\ 
\par
Besides, we get: $\psi ^{m+1}(u_{n},v_{n})=(m+1)\Phi _{n},\ n\geq 0.$\ 
\par
The estimate of the numerical blow-up is in that case the same as the 
estimate of the exact blow-up.\ \par
\ \par
\subsection{Properties of the scheme}
\ \par
\begin{Proposition} If$\ \alpha >\lambda _{1}$ and  $T^{*}$ is the numerical 
blow-up time, we get the estimate: \ \par

\begin{displaymath} 
\Phi _{n}^{\frac{1}{m+1}}\leq \left({\frac{T_{*}}{T_{*}-t}}\right) ^{\frac{1}{1-m}}\Phi 
_{0}^{\frac{1}{m+1}}\end{displaymath} \ \par
\end{Proposition}
{\it{Proof:}} The estimate ~(\ref{systeme26}) may be written in this 
case:\ \par
$(1-m)\Delta tF(u_{n+1},v_{n+1})\leq \psi ^{m-1}(u_{n+1},v_{n+1})-\psi 
^{m-1}(u_{n},v_{n})\leq (1-m)\Delta tF(u_{n},v_{n}).$\ \par
\ \par
We deduce, since $(F(u_{n},v_{n}))_{n\geq 0}$ is a nonincreasing sequence:\ 
\par
\ \par
$\psi _{n}^{m-1}(u_{n+j},v_{n+j})-\psi _{n}^{m-1}(u_{n},v_{n})\leq (1-m)j\Delta 
tF(u_{n},v_{n}),\ \forall j\geq 0,$\ \par
\ \par
$\psi _{n}^{m-1}(u_{n-i},v_{n-i})-\psi _{n}^{m-1}(u_{n},v_{n})\leq -(1-m)i\Delta 
tF(u_{n},v_{n}),\ \forall i\geq 0.$\ \par
So, we obtain: \ \par
\ \par
$i\psi _{n}^{m-1}(u_{n+j},v_{n+j})+j\psi _{n}^{m-1}(u_{n-i},v_{n-i})\leq 
(i+j)\psi _{n}^{m-1}(u_{n},v_{n}).$\ \par
\ \par
If $T_{*}=N\Delta t$ is the numerical blow-up time, by choosing $i=n,\ j=N-n$, 
we get :\ \par
$(N-n)\psi ^{m-1}(u_{0},v_{0})\leq N\psi ^{1-m}(u_{n},v_{n})$ or, $\displaystyle 
\psi ^{1-m}(u_{n},v_{n})\leq \frac{\displaystyle T_{*}}{\displaystyle 
T_{*}-t}\psi ^{1-m}(u_{n},v_{n})$.\ \par
This is the same estimate as for the exact solution.\ \par
\ \par
\ \par
\begin{Proposition} If $\alpha <\lambda _{1}$, the numerical solution 
tends to $0\ \ $when $t{\longrightarrow}+\infty $.\ \par
\end{Proposition}
{\it{Proof: }}In that case, we get:{\it{ }}$J(u,v)\geq (\lambda _{1}-\alpha )\int_{\Omega }^{}{(u^{2}+v^{2})dx}\geq 
0.$\ \par
Besides, since $m<1$, for any $\phi \in L^{2}(\Omega )$, we have: $\int_{\Omega }^{}{\phi ^{2}dx}\geq C(\Omega )\left({\int_{\Omega }^{}{\phi 
^{m+1}dx}}\right) ^{\frac{2}{m+1}}$\ \par
and by using Young inequality, we get \ \par
$\displaystyle 2^{\displaystyle -\frac{\displaystyle 1-m}{\displaystyle m+1}}\displaystyle 
\left({\displaystyle \displaystyle \int_{\Omega }^{}{(u^{m+1}+v^{m+1})dx}}\right) 
^{\displaystyle ^{\displaystyle \frac{\displaystyle 2}{\displaystyle 
m+1}}}\leq \displaystyle \left({\displaystyle \displaystyle \int_{\Omega 
}^{}{u^{m+1}dx}}\right) ^{\displaystyle \frac{\displaystyle 2}{\displaystyle 
m+1}}+\displaystyle \left({\displaystyle \displaystyle \int_{\Omega }^{}{v^{m+1}dx}}\right) 
^{\displaystyle \frac{\displaystyle 2}{\displaystyle m+1}}$\ \par
$\displaystyle \leq \frac{\displaystyle 1}{\displaystyle C(\Omega )}\displaystyle \left({\displaystyle 
\displaystyle \int_{\Omega }^{}{(u^{2}+v^{2})dx}}\right) $.\ \par
So, we get : $\displaystyle F(u,v)\geq \frac{\displaystyle \lambda _{1}-\alpha }{\displaystyle C(\Omega 
)}.$\ \par
Besides, from ~(\ref{systeme26}), we obtain $\displaystyle \psi ^{m-1}(u_{n+1},v_{n+1})\geq \psi ^{m-1}(u_{n},v_{n})\geq \frac{\displaystyle 
\lambda _{1}-\alpha }{\displaystyle C(\Omega )}\Delta t$ and\ \par
$\displaystyle \psi ^{m-1}(u_{n},v_{n})\geq \psi ^{m-1}(u_{0},v_{0})+\frac{\displaystyle 
\lambda _{1}-\alpha }{\displaystyle C(\Omega }t_{n}$.\ \par
We deduce: $\underset{n{\longrightarrow}+\infty }{\lim }\int_{\Omega }^{}{(u^{m+1}+v^{m+1})dx}=0$. 
\ \par
\ \par
\begin{Proposition} If $\alpha =\lambda _{1},\ then\ \underset{n{\longrightarrow}+\infty }{\lim 
}u_{n}=\underset{n{\longrightarrow}+\infty }{\lim }v_{n}=\theta \rho 
_{1}$ where $\theta $ is depending on the initial conditions.\ \par
\end{Proposition}
{\it{Proof:}} We have the equality: $\Phi _{n+1}-\Phi _{n}=\frac{1}{m+1}\left({\psi ^{m+1}(u_{n+1},v_{n+1})-\psi 
^{m+1}(u_{n},v_{n})}\right) $\ \par
and from the inequality $\displaystyle \forall a,b\in {\mathbb{R}}^{+},\ a^{m+1}-b^{m+1}\leq \frac{\displaystyle 
1+m}{\displaystyle 1-m}a^{2}\displaystyle \left({\displaystyle b^{m-1}-a^{m-1}}\right) 
$, \ \par
we obtain:\ \par
$\displaystyle \Phi _{n+1}-\Phi _{n}\leq \frac{\displaystyle 1}{\displaystyle 1-m}\psi 
^{2}(u_{n+1},v_{n+1})\displaystyle \left({\displaystyle \psi ^{m-1}(u_{n},v_{n})-\psi 
^{m-1}(u_{n+1},v_{n+1})}\right) \leq -\Delta tJ(u_{n+1},v_{n+1})$.\ \par
\ \par
In the same manner as in proposition ~(\ref{systeme33}), we deduce that 
$\underset{n{\longrightarrow}+\infty }{\lim }u_{n}=\underset{n{\longrightarrow}+\infty 
}{\lim }v_{n}=\theta \rho _{1}$.\ \par
We now obtain an estimate for $\theta :$\ \par
By multiplying the equations ~(\ref{systeme11}), ~(\ref{systeme12}) by 
$\rho _{1}$ and integrating on $\Omega ,$ we get:\ \par
\ \par
$\displaystyle \frac{\displaystyle m}{\displaystyle 1-m}\displaystyle \int_{\Omega }^{}{u_{n+1}u_{n}^{m-1}\rho 
_{1}dx}-\frac{\displaystyle m}{\displaystyle 1-m}\displaystyle \int_{\Omega 
}^{}{u_{n+1}^{m}\rho _{1}dx}+\lambda _{1}\Delta t\displaystyle \int_{\Omega 
}^{}{(u_{n+1}-v_{n+1})\rho _{1}dx}=0,$\ \par
\ \par
$\displaystyle \frac{\displaystyle m}{\displaystyle 1-m}\displaystyle \int_{\Omega }^{}{v_{n+1}v_{n}^{m-1}\rho 
_{1}dx}-\frac{\displaystyle m}{\displaystyle 1-m}\displaystyle \int_{\Omega 
}^{}{v_{n+1}^{m}\rho _{1}dx}+\lambda _{1}\Delta t\displaystyle \int_{\Omega 
}^{}{(v_{n+1}-u_{n+1})dx}=0.$\ \par
\ \par
By using the inequality: $u_{n+1}^{m}\leq (1-m)u_{n}^{m}+mu_{n+1}u_{n}^{m-1}$ 
and the same inequality applied to the function$\ v,$ we get \ \par
\ \par
$\displaystyle \displaystyle \int_{\Omega }^{}{u_{n+1}^{m}\rho _{1}dx}\leq \displaystyle 
\int_{\Omega }^{}{u_{n}^{m}\rho _{1}dx}+\Delta t\frac{\displaystyle \lambda 
_{1}}{\displaystyle m}\displaystyle \int_{\Omega }^{}{(v_{n+1}-u_{n+1})\rho 
_{1}dx},$\ \par
\ \par
$\displaystyle \displaystyle \int_{\Omega }^{}{v_{n+1}^{m}\rho _{1}dx}\leq \displaystyle 
\int_{\Omega }^{}{v_{n}^{m}\rho _{1}dx}+\Delta t\frac{\displaystyle \lambda 
_{1}}{\displaystyle m}\displaystyle \int_{\Omega }^{}{(u_{n+1}-v_{n+1})\rho 
_{1}dx},$\ \par
and then: $\displaystyle \displaystyle \int_{\Omega }^{}{(u_{n+1}^{m}+v_{n+1}^{m})\rho _{1}dx}\leq 
\displaystyle \int_{\Omega }^{}{(u_{0}^{m}+v_{0}^{m})\rho _{1}dx},\ n\geq 
0$.\ \par
\ \par
We deduce: $\displaystyle \theta ^{m}\leq \frac{\displaystyle \displaystyle \int_{\Omega }^{}{(u_{0}^{m}+v_{0}^{m})\rho 
_{1}dx}}{\displaystyle 2\displaystyle \int_{\Omega }^{}{\rho _{1}^{m+1}dx}}$.\ 
\par
\ \par
\subsection{Convergence of the scheme}
\ \par
In this section, we obtain estimates on the numerical solution so we 
can extract by compactness a convergent subsequence. In order to prove 
that the limit is solution of the system ~(\ref{systeme0}), we need an 
hypothesis on the initial condition, (this is due to the fact that we 
have a negative power in the scheme and we may not use a Hold{\"e}r inequality). 
If the initial condition does not satisfy the hypothesis, we observe 
numerically that this hypothesis is satisfied after a few times steps 
and the scheme again converges.\ \par
\ \par
Let us denote by $T_{1}^{*}=\underset{0<\Delta t<\Delta t_{0}}{\inf }T_{*}(\Delta t)$ 
if $T_{*}$ is the existence time of the numerical solution. It follows 
from theorem ~\ref{systeme34} that $T_{1}^{*}\geq T_{1}>0.$ Let $T\in [0,T_{1}^{*}[,\ (T=N\Delta t).$ 
We denote $u_{\Delta t}$ and $v_{\Delta t}$ the approximation of $u$ 
and $v$ defined by:\ \par

\begin{displaymath} 
u_{\Delta t}(t)=\left({u_{n}^{m}+\frac{t-t_{n}}{\Delta t}(u_{n+1}^{m}-u_{n}^{m})}\right) 
^{1/m},\end{displaymath} \ \par

\begin{displaymath} 
v_{\Delta t}(t)=\left({v_{n}^{m}+\frac{t-t_{n}}{\Delta t}(v_{n+1}^{m}-v_{n}^{m})}\right) 
^{1/m}.\end{displaymath} \ \par
\ \par
\begin{Theorem} The sequences $(u_{\Delta t})$ and $(v_{\Delta t})$ are 
uniformly bounded in $C(0,T;H_{0}^{1}(\Omega ))$ and $H^{1}(0,T;L^{2}(\Omega )).$\ 
\par
\end{Theorem}
\ \par
\textsl{Proof: } Since $T<T_{1}^{*}$, the functions $(u_{n})_{n\geq 0}$ 
 and $(v_{n})_{n\geq 0}$  are uniformly bounded in $C(0,T;\Omega )$ and 
since $J(u_{n},v_{n})$ is nonincreasing, we get: $J(u_{n},v_{n})\leq J(u_{0},v_{0}).$ 
We deduce that the sequences $(\nabla u_{n})_{n\geq 0}$  and $(\nabla v_{n})_{n\geq 0}$ 
are uniformly bounded in $L^{2}(\Omega ).$\ \par
We prove now that the sequences $\displaystyle (\frac{\displaystyle du_{\Delta t}}{\displaystyle dt})$ 
and $(\frac{dv_{\Delta t}}{dt})$ are uniformly bounded in $L^{2}(0,T;L^{2}(\Omega )).$\ 
\par
We have: \ \par

\begin{displaymath} 
\displaystyle \left\Vert{\displaystyle \frac{\displaystyle du_{\Delta 
t}}{\displaystyle dt}}\right\Vert ^{2}_{\displaystyle L^{2}(0,T;L^{2}(\Omega 
))}=\displaystyle \sum_{n=0}^{N-1}{\displaystyle \int_{t_{n}}^{t_{n+1}}{\displaystyle 
\int_{\Omega }^{}{\frac{\displaystyle 1}{\displaystyle m^{2}}\ u_{\Delta 
t}^{2(1-m)}\displaystyle \left({\displaystyle \frac{\displaystyle u_{n+1}^{m}-u_{n}^{m}}{\displaystyle 
\Delta t}}\right) ^{2}dx}dt}}\end{displaymath} \ \par
\ \par

\begin{displaymath} 
\leq \frac{\displaystyle 1}{\displaystyle m^{2}\Delta t}\displaystyle 
\sum_{n=0}^{N-1}{\displaystyle \int_{\Omega }^{}{\displaystyle \left({\displaystyle 
u_{n+1}^{2(1-m)}+u_{n}^{2(1-m)}}\right) }}(u_{n+1}^{m}-u_{n}^{m})^{2}dx.\end{displaymath} 
\ \par
\ \par
\ \par
From the following inequalities, $u_{n+1}^{1-m}\left\vert{u_{n+1}^{m}-u_{n}^{m}}\right\vert \leq \left\vert{u_{n+1}-u_{n}}\right\vert 
$ and $u_{n}^{1-m}\left\vert{u_{n+1}^{m}-u_{n}^{m}}\right\vert \leq \left\vert{u_{n+1}-u_{n}}\right\vert 
$, we deduce:\ \par

\begin{displaymath} 
\displaystyle \left\Vert{\displaystyle \frac{\displaystyle du_{\Delta 
t}}{\displaystyle dt}}\right\Vert ^{2}_{\displaystyle L^{2}(0,T;L^{2}(\Omega 
))}\leq \frac{\displaystyle 1}{\displaystyle m^{2}\Delta t}\displaystyle 
\sum_{n=0}^{N-1}{\displaystyle \int_{\Omega }^{}{(u_{n+1}-u_{n})^{2}dx}}\end{displaymath} 
\ \par
\ \par

\begin{displaymath} 
\leq \frac{\displaystyle 1}{\displaystyle m^{2}\Delta t}\displaystyle 
\sum_{n=0}^{N-1}{\displaystyle \left({\displaystyle \displaystyle \int_{\Omega 
}^{}{u_{n}^{m-1}(u_{n+1}-u_{n})^{2}dx}}\right) \displaystyle \left\Vert{\displaystyle 
u_{n}}\right\Vert _{\infty }^{1-m}}\end{displaymath} \ \par
and we obtain analogous inequality for $v_{\Delta t}$.\ \par
Hence, we get: \ \par

\begin{displaymath} 
\left\Vert{\frac{du_{\Delta t}}{dt}}\right\Vert ^{2}_{L^{2}(0,T;L^{2}(\Omega 
)}+\left\Vert{\frac{dv_{\Delta t}}{dt}}\right\Vert ^{2}_{L^{2}(0,T;L^{2}(\Omega 
))}\leq C\ (J(u_{0},v_{0})-J(u_{N},v_{N}))\end{displaymath} \ \par
since $(\left\Vert{u_{n}}\right\Vert _{\infty })_{0\leq n\leq N\ }$ and 
$(\left\Vert{v_{n}}\right\Vert _{\infty })_{0\leq n\leq N}\ $ are uniformly 
bounded.\ \par
If $\alpha \leq \lambda _{1},$ we have: $J(u_{N},v_{N})\geq 0,$\ \par
if $\alpha >\lambda _{1},$ we get $J(u_{N},v_{N})\geq (\lambda _{1}-\alpha )\int_{\Omega }^{}{(u_{N}^{2}+v_{N}^{2})dx}\ 
$ and this quantity is bounded from below; we deduce that the sequences 
are uniformly bounded in $L^{2}(0,T;L^{2}(\Omega )).$\ \par
\ \par
Since the sequences $(u_{\Delta t})$ and $(v_{\Delta t})$ are uniformly 
bounded in $C(0,T;H_{0}^{1}(\Omega ))\cap H^{1}(0,T;L^{2}(\Omega )),$ 
we can extract subsequences which converge to functions $u$ and $v$ in 
$C(0,T;\Omega )$ if $d=1$ and in $C(0,T;L^{r}(\Omega )),\ (r<2d/(d-2)$ 
if $d>2$, $r=\infty $ if $d=2,$( Simon ~\cite{Simon}).\ \par
In order to prove that the limits $u,\ v$ are solutions of the system 
~(\ref{systeme0}), we use the same proof as in ~\cite{mnLeRoux1}; so 
we need to estimate the quantities $\left\vert{(u_{n+1}^{1-m}-u_{n}^{1-m})u_{n}^{m-1}}\right\vert $ 
and $\left\vert{(v_{n+1}^{1-m}-v_{n}^{1-m})v_{n}^{m-1}}\right\vert $. 
This is the object of the two next lemmas.\ \par
\ \par
\begin{Lemma} For $n\geq 0,$ we have the inequalities:\ \par

\begin{equation} 
t_{n+1}u_{n+1}^{1-m}\geq t_{n}u_{n}^{1-m},\ t_{n+1}v_{n+1}^{1-m}\geq t_{n}v_{n}^{1-m}.\label{systeme35}
\end{equation} \ \par
\end{Lemma}
\textsl{Proof:} These inequalities are proved recurrently; it is true 
for $n=0.$ Assume it is true at the order $n-1,$ that is :\ \par
$t_{n}u_{n}^{1-m}\geq t_{n-1}u_{n-1}^{1-m},\ t_{n}v_{n}^{1-m}\geq t_{n-1}v_{n-1}^{1-m}.$\ 
\par
The functions $\left({\frac{t_{n}}{t_{n+1}}}\right) ^{1/(1-m)}u_{n}$ 
, $\left({\frac{t_{n}}{t_{n+1}}}\right) ^{1/(1-m)}v_{n}$ will be subsolutions 
of ~(\ref{systeme0}) if :\ \par

\begin{displaymath} 
-\frac{m}{1-m}\frac{\Delta t}{t_{n}}u_{n}^{m}+\Delta tAu_{n}-\alpha \Delta 
tv_{n}\leq 0,\end{displaymath} \ \par

\begin{displaymath} 
-\frac{m}{1-m}\frac{\Delta t}{t_{n}}v_{n}^{m}+\Delta tAv_{n}-\alpha \Delta 
tu_{n}\leq 0.\end{displaymath} \ \par
But $u_{n},v_{n}$ are the solutions at the time level $t_{n}$, so we 
get:\ \par
$\displaystyle -\frac{\displaystyle m}{\displaystyle 1-m}\frac{\displaystyle \Delta t}{\displaystyle 
t_{n}}u_{n}^{m}+\Delta t(Au_{n}-\alpha v_{n})=-\frac{\displaystyle m}{\displaystyle 
1-m}\displaystyle \left({\displaystyle \frac{\displaystyle \Delta t}{\displaystyle 
t_{n}}u_{n}^{m}+u_{n}u_{n-1}^{m-1}-u_{n}^{m}}\right) =-\frac{\displaystyle 
m}{\displaystyle 1-m}u_{n}^{m}\displaystyle \left({\displaystyle u_{n}^{1-m}u_{n-1}^{m-1}-\frac{\displaystyle 
t_{n-1}}{\displaystyle t_{n}}}\right) $\ \par
and $\displaystyle -\frac{\displaystyle m}{\displaystyle 1-m}\frac{\displaystyle \Delta t}{\displaystyle 
t_{n}}v_{n}^{m}+\Delta t(Av_{n}-\alpha u_{n})=-\frac{\displaystyle m}{\displaystyle 
1-m}v_{n}^{m}(v_{n}^{1-m}v_{n-1}^{m-1}-\frac{\displaystyle t_{n-1}}{\displaystyle 
t_{n}}).$\ \par
From the recurrence hypotheses, the second members of these two inequalities 
are non positive and the inequalities~(\ref{systeme35}) are satisfied 
at the order $n.$\ \par
\ \par
\begin{Lemma} Assume that the initial conditions $u_{0},v_{0}$ satisfy:\ 
\par
There exists a constant $C_{0}$ such that \ \par

\begin{equation} 
Au_{0}-\alpha v_{0}+C_{0}u_{0}^{m}\geq 0,\ Av_{0}-\alpha u_{0}+C_{0}v_{0}^{m}\geq 
0\label{systeme36}
\end{equation} \ \par
Then, we have the estimates: \ \par

\begin{equation} 
u_{n+1}\leq \left({\frac{T_{2}-t_{n}}{T_{2}-t_{n+1}}}\right) ^{1/(1-m)}u_{n},\ 
v_{n+1}\leq \left({\frac{T_{2}-t_{n}}{T_{2}-t_{n+1}}}\right) ^{1/(1-m)}v_{n}\label{systeme37}
\end{equation} \ \par
with $T_{2}=\frac{m}{(1-m)C_{0}}$\ \par
\end{Lemma}
\ \par
\textsl{Proof: }This lemma is proved recursively. First, we prove the 
inequalities for $n=0.$\ \par
$\left({\frac{T_{2}}{T_{2}-\Delta t}}\right) ^{1/(1-m)}u_{0},\ \left({\frac{T_{2}}{T_{2}-\Delta 
t}}\right) ^{1/(1-m)}\ $will be supersolutions of ~(\ref{systeme11}), 
~(\ref{systeme12}) if \ \par
\ \par

\begin{displaymath} 
-\frac{m}{1-m}\left({\frac{T_{2}}{T_{2}-\Delta t}}\right) ^{m/(1-m)}u_{0}^{m}+\frac{m}{1-m}\left({\frac{T_{2}}{T_{2}-\Delta 
t}}\right) ^{1/(1-m)}u_{0}^{m}+\Delta t\left({\frac{T_{2}}{T_{2}-\Delta 
t}}\right) ^{1/(1-m)}(Au_{0}-\alpha v_{0})\geq 0,\end{displaymath} \ 
\par
\ \par

\begin{displaymath} 
-\frac{m}{1-m}\left({\frac{T_{2}}{T_{2}-\Delta t}}\right) ^{m/(1-m)}v_{0}^{m}+\frac{m}{1-m}\left({\frac{T_{2}}{T_{2}-\Delta 
t}}\right) ^{1/(1-m)}v_{0}^{m}+\Delta t\left({\frac{T_{2}}{T_{2}-\Delta 
t}}\right) ^{1/(1-m)}(Av_{0}-\alpha u_{0})\geq 0,\end{displaymath} \ 
\par
\ \par
that is,    $\frac{m}{(1-m)T_{2}}u_{0}^{m}+Au_{0}-\alpha v_{0}\geq 0$ 
and    $\frac{m}{(1-m)T_{2}}v_{0}^{m}+Av_{0}-\alpha u_{0}\geq 0$.  \ 
\par
\ \par
From ~(\ref{systeme36}), these inequalities are satisfied.\ \par
Assume now that the estimates are satisfied at the order $n.$We prove 
it is satisfied at the order $n+1.$ \ \par
The functions $\left({\frac{T_{2}-t_{n}}{T_{2}-t_{n+1}}}\right) ^{1/(1-m)}u_{n},\ \left({\frac{T_{2}-t_{n}}{T_{2}-t_{n+1}}}\right) 
^{1/(1-m)}v_{n}$  will be supersolutions of ~(\ref{systeme11}), ~(\ref{systeme12}) 
if \ \par
$\displaystyle \frac{\displaystyle m}{\displaystyle 1-m}\frac{\displaystyle 1}{\displaystyle 
T_{2}-t_{n}}u_{n}^{m}+Au_{n}-\alpha v_{n}\geq 0,\ \frac{\displaystyle 
m}{\displaystyle 1-m}\frac{\displaystyle 1}{\displaystyle T_{2}-t_{n}}v_{n}^{m}+Av_{n}-\alpha 
u_{n}\geq 0.$\ \par
Since $u_{n},\ v_{n}$ are the solutions at the time level $t_{n}$, we 
get :\ \par

\begin{displaymath} 
\frac{m}{1-m}\frac{u_{n}^{m}}{T_{2}-t_{n}}+Au_{n}-\alpha v_{n}=\frac{m}{(1-m)\Delta 
t}\left({u_{n}^{m}\frac{T_{2}-t_{n-1}}{T_{2}-t_{n}}-u_{n}u_{n-1}^{m-1}}\right) 
,\end{displaymath} \ \par

\begin{displaymath} 
\frac{m}{1-m}\frac{v_{n}^{m}}{T_{2}-t_{n}}+Av_{n}-\alpha u_{n}=\frac{m}{(1-m)\Delta 
t}\left({v_{n}^{m}\frac{T_{2}-t_{n-1}}{T_{2}-t_{n}}-v_{n}v_{n-1}^{m-1}}\right) 
.\end{displaymath} \ \par
\ \par
By using the recurrence hypothesis, we obtain the estimates ~(\ref{systeme37}).\ 
\par
\ \par
\begin{Remark}  If $Au_{0}-\alpha v_{0}\geq 0$ and $Av_{0}-\alpha u_{0}\geq 0$, 
then $T_{2}=+\infty \ $and we obtain: $u_{n+1}\leq u_{n},\ v_{n+1}\leq v_{n}.$\ 
\par
\end{Remark}
\ \par
From these two lemmas, we obtain the inequalities \ \par

\begin{displaymath} 
-\frac{\Delta t}{t_{n+1}}\leq (u_{n+1}^{1-m}-u_{n}^{1-m})u_{n}^{m-1}\leq 
\frac{\Delta t}{T_{2}-t_{n+1}}\end{displaymath} \ \par
and analogous inequalities for $v_{n}.$ Then, we obtain the convergence 
of the scheme as in ~\cite{mnLeRoux1}.\ \par

\bibliographystyle{C:/TexLive/texmf/bibtex/bst/base/plain}

\end{document}